\documentclass[12pt,a4paper]{article}

\usepackage{amsmath}    
\usepackage{amssymb}
\usepackage{graphicx}   
\usepackage{verbatim}   
\usepackage{color}      
\usepackage{hyperref}   
\usepackage{enumerate}
\usepackage{listings}
\usepackage{xcolor}
\usepackage{slashed}
\usepackage{amsfonts}
\usepackage{amsmath}    
\usepackage{amssymb}
\usepackage{amsthm}
\usepackage{mathrsfs}
\usepackage{empheq}
\usepackage{bbold}
\usepackage{etoolbox}

\usepackage{cleveref}
\crefformat{footnote}{#2\footnotemark[#1]#3}

\usepackage[british]{babel}
\usepackage[font={small,it}]{caption}

\newcommand{\del}{\partial}
\renewcommand{\theta}{\vartheta}
\renewcommand{\phi}{\varphi}

\newcommand{\dd}{\mathrm{d}}

\newcommand{\id}{\mathbb{1}}

\renewcommand{\vec}{\mathbf}

\renewcommand{\chi}{\id}

\definecolor{mygreen}{rgb}{0,0.6,0}
\definecolor{mygray}{rgb}{0.5,0.5,0.5}
\definecolor{mymauve}{rgb}{0.58,0,0.82}

\newtheorem{definition}{Definition}[section]

\newtheorem{theorem}{Theorem}[section]
\newtheorem{example}{Example}[section]

\AtEndEnvironment{example}{%
  \qed
}

\renewcommand{\title}{Semi-discrete Active Flux as a Petrov-Galerkin method: \\the case of one-dimensional and Cartesian grids}

\newcommand{\authorOne}{Wasilij Barsukow
\footnote{Bordeaux Institute of Mathematics, Bordeaux University and CNRS/UMR5251, Talence, 33405 France, wasilij.barsukow@math.u-bordeaux.fr}
}
\begin{document}

\begin{center} \Large
\title

\vspace{1cm}

\date{}
\normalsize

\authorOne
\end{center}

\begin{abstract}

% #############################################################################################################################
Active Flux (AF) is a numerical method for hyperbolic conservation laws, whose degrees of freedom are averages/moments and (shared) point values at cell interfaces. It has been noted previously in a heuristic fashion that it thus combines ideas from Finite Volume/Discontinuous Galerkin (DG) methods with a continuous approximation common in continuous Finite Element (CG) methods. This work shows that semi-discrete Active Flux methods can be obtained from a variational formulation through a particular choice of (biorthogonal) test functions. These latter being discontinuous, the new formulation emphasizes the intermediate nature of AF between DG and CG. Explicit constructions are given for the case of arbitrarily high-order Active Flux with additional moments in 1-d, and for the classical third-order Active Flux on two-dimensional Cartesian meshes.
% #############################################################################################################################



\end{abstract}

\section{Introduction}

Active Flux is a numerical method for conservation laws in one
\begin{align}
\del_t q + \del_x f(q) &= 0 & q &\colon \mathbb R^+_0 \times \mathbb R \to \mathbb R^m \label{eq:conslaw1d}\\
&&f &\colon \mathbb R^m \to \mathbb R^m \nonumber
\end{align}
and several
\begin{align}
\del_t q + \del_x f(q) + \del_y g(q) &= 0 & q&\colon \mathbb R^+_0 \times \mathbb R^2 \to \mathbb R^m \label{eq:conslaw2d}\\
&&f,g &\colon \mathbb R^m \to \mathbb R^m \nonumber
\end{align}
space dimensions.
In its third-order version (\cite{vanleer77,eymann13}), the method makes use of cell averages and point values at the boundaries of the computational cell as independent degrees of freedom. In \cite{abgrall22,lechner25} arbitrarily high-order extensions were proposed with moments as further degrees of freedom. Active Flux approximates the solution in a globally continuous way, since the point values at cell interfaces are shared between adjacent cells. Riemann problems therefore do not arise at cell interfaces and the flux can be evaluated directly. The degrees of freedom are chosen such that in every cell a polynomial of high degree can be reconstructed, while the global continuity is merely $C^0$.

So far, in the literature, Active Flux has not been presented as originating in a variational formulation.
The averages can be updated immediately by evaluating the flux using the point values at cell interfaces.
The traditional way of updating the point values assumes that the initial value problem for the given PDE with piecewise polynomial, continuous initial data can be solved exactly or sufficiently accurately. Exact evolution operators are only available for linear problems such as linear advection (\cite{vanleer77}) and linear acoustics (\cite{fan15,barsukow18activeflux}). For nonlinear scalar conservation laws and for hyperbolic systems in one spatial dimension, approximate evolution operators are known (\cite{kerkmann18,barsukow19activeflux}, and recently there has been some progress on evolution operators for multi-dimensional nonlinear systems (\cite{chudzik24,barsukow25afos,chudzik25}). 

In \cite{abgrall20} it has been proposed to use the degrees of freedom to discretize the PDE in space first, and then apply a method-of-line approach to the resulting ODEs. Again, the update equation for the averages is obvious (see below: Equations \eqref{eq:averagesAF} for 1-d, and \eqref{eq:averagesAF2d} for 2-d) and in multi-d only requires choosing a suitable quadrature along the edge. For the point value, the space derivative is replaced by a finite difference formula. A straightforward choice is to evaluate the derivative of the reconstruction at the location of the point value. Since, globally, the reconstruction is only $C^0$, upwinding is included by differentiating the reconstruction in the upwind cell. A central update is obtained by taking the average of the derivatives from the two adjacent cells. In \cite{abgrall24pampa} it has been shown how the methodology of Virtual Finite Elements can be used to obtain reconstructions on general polygons.

The semi-discrete Active Flux has been shown to be stationarity and vorticity preserving for linear acoustics, and to be able to cope both with strong shocks and the low Mach number regime (\cite{duan24,barsukow24afeuler}). Shared degrees of freedom make it significantly more efficient than Discontinuous Galerkin (DG) methods since there are less values to update and since point values are cheaper to update than moments. At the same time, all operations remain local. A performance comparison is presented in \cite{barsukow26afdg}. It thus seems very well suited for efficiently solving multi-dimensional conservation laws such as the Euler or MHD equations (see e.g. \cite{duan25mhd,liu25mhd}).

While the arbitrary-order extension using higher moments bears a certain similitude to DG methods (method-of-lines, integration-by-parts), the following are significant differences:
\begin{itemize}
 \item It is mass-matrix-free\footnote{By ``mass-matrix-free'', or ``absence of mass matrix'' we mean here and in the following that the mass matrix is the identity matrix.}, with not even a small local mass matrix to invert.
 \item The solution is globally continuous, i.e. the point values at cell interfaces are shared among the adjacent cells.
\end{itemize}
Although not a conceptual difference, the fact that both nodal and modal degrees of freedom are simultaneously present is a non-classical choice in DG, too. Global continuity makes Active Flux be somewhat close to continuous Finite Element methods, but other of its features are again untypical, such as the absence of a variational formulation.

This paper aims at showing that the semi-discrete Active Flux can be very naturally obtained as a Petrov-Galerkin method. The continuous basis functions are very similar to what is customary for Continuous Galerkin (CG) methods, while the test functions will be chosen such that there is no mass matrix and will turn out to be discontinuous at cell interfaces, as are those of DG. The specific examples considered here are one-dimensional Active Flux methods of arbitrarily high order (via additional moments) and the classical third-order Active Flux method on two-dimensional Cartesian grids. Thus, in the end the antagonism is resolved since these Active Flux methods combine the different features of CG and DG in a non-contradicting way. It is shown that upwinding comes about because of discontinuities in the test functions.

A mass matrix which is the identity is obtained by choosing so-called biorthogonal test functions. Recall that families of orthogonal polynomials are usually constructed such that the polynomial $p_n$ of some degree $n$ is orthogonal (for some scalar product $(\,\cdot \,, \, \cdot \,)$) to all polynomials of lesser degree: $(p_n, p_m) = 0$ $\forall m < n$. Biorthogonal polynomial families (as studied e.g. in \cite{konhauser65}) are two families $(p_n)_n, (q_n)_n$ such that $(p_n, q_m) = 0 $ $\forall m \neq n$. The natural ordering of the two families by degree of the polynomial, however, is not what is needed in the context of finite elements: e.g. nodal basis functions typically are all of the same degree. Results on biorthogonal polynomial families, therefore, are of restricted use here. However, as is discussed below, given a set of basis functions $(\phi_k)_k$, it can be possible to find test functions $(\psi_k)_k$ such that $(\psi_k, \phi_\ell) = \delta_{k\ell}$, i.e. such that the mass matrix is the identity. Such biorthogonal test functions have been used e.g. in \cite{wohlmuth00}; see Section \ref{sec:biorthogonal} for further references.

The paper is organized as follows. Section \ref{sec:1d} gives an overview of Active Flux and introduces the ideas on biorthogonal test and basis functions in one spatial dimension. It shows that the Petrov-Galerkin approach leads directly to the Active Flux method. Section \ref{sec:2d} demonstrates the same approach on 2-d Cartesian grids. The 2-d basis functions are not a tensor product of one-dimensional basis functions, which makes an explicit analysis necessary. At the end of the Section, the relation to a tensor basis is elucidated. 

In one spatial dimension, the computational grid is assumed to consist of cells $C_i = [x_{i-\frac12},x_{i+\frac12})$, $i \in \mathbb Z$ in the one-dimensional case and of cells $C_{ij} = [x_{i-\frac12},x_{i+\frac12}) \times [y_{j-\frac12},y_{j+\frac12})$, $(i,j) \in \mathbb Z^2$ in two dimensions. The grids are assumed equidistant with $\Delta x := x_{i+\frac12}-x_{i-\frac12}$ and $\Delta y := y_{j+\frac12}-y_{j-\frac12}$. Denote by $P^K$ the space of univariate polynomials of degree at most $K$, and by $Q^K \equiv P^{K,K}$ the tensor-product space of polynomials in $x$ and $y$ obtained from $P^K$. The $L^2$ scalar product is denoted by
\begin{align}
 (u,v)_{L^2} := \int u(x) v(x) \, \dd x
\end{align}
in 1-d and by
\begin{align}
 (u,v)_{L^2} := \iint u(x,y) v(x,y) \, \dd x\dd y
\end{align}
in 2-d.
$\chi_I$ is the characteristic function of the interval $I$.

\section{One spatial dimension} \label{sec:1d}

\newcommand{\mombas}{A}

\subsection{Review of the semi-discrete Active Flux method}

This Section gives a brief overview of the key properties of the Active Flux method in one space dimension.
In its simplest form, the Active Flux method uses cell averages and point values at cell interfaces
\begin{align}
 \bar q_i(t) &\simeq \frac{1}{\Delta x} \int_{x_{i-\frac12}}^{x_{i+\frac12}} q(t, x) \,\dd x & q_{i+\frac12}(t) \simeq q(t, x_{i+\frac12})
\end{align}
as degrees of freedom. Upon integration of the conservation law \eqref{eq:conslaw1d} over the cell and using Gauss' law one immediately finds the evolution equation for the cell average:
\begin{align}
 \frac{\dd}{\dd t} \bar q_i + \frac{f(q_{i+\frac12}) - f(q_{i-\frac12})}{\Delta x} &= 0 \label{eq:averagesAF}
\end{align}
Observe that the computation of the numerical flux does not require a Riemann solver, since data is directly available at the cell interface.

In \cite{abgrall20,abgrall22} it has been suggested to update the point values by first replacing $\del_x f(q)$ by $J(q) \del_x q$ with the Jacobian $J = \nabla_q f$, and then by further splitting the expression into left-going and right-going waves:
\begin{align}
 \frac{\dd}{\dd t} q_{i+\frac12} + J^+ (D^+ q)_{i+\frac12} +  J^- (D^- q)_{i+\frac12} &= 0 \label{eq:pointvaluesAF}
\end{align}
Here, $J^\pm := \frac12 (J(q_{i+\frac12}) \pm \vert J(q_{i+\frac12})\vert )$ (where the absolute value is computed on the eigenvalues), and $D^{\pm} q$ are finite differences biased to the left/right, respectively. If $(D^+ q)_{i+\frac12}$ depends only on $q_{i-\frac12}, \bar q_i, q_{i+\frac12}$ (and $(D^- q)_{i+\frac12}$ only on $q_{i+\frac12}, \bar q_{i+1}, q_{i+\frac32}$), then they are unique and are equal to the derivatives of quadratic interpolations of the three values involved. One is led to consider in every cell $i$ the unique quadratic polynomial $q_{\text{recon},i} \colon \left[ -\frac{\Delta x}{2}, \frac{\Delta x}{2} \right] \to \mathbb R$ which satisfies
\begin{align}
 \frac{1}{\Delta x} \int_{x_{i-\frac12}}^{x_{i+\frac12}}q_{\text{recon},i} \,\dd x  &= \bar q_i & q_{\text{recon},i}\left( \pm \frac{\Delta x}{2} \right) &= q_{i\pm\frac12}
\end{align}

These polynomials, as is customary for Finite Volume methods, are called reconstructions in existing literature. The derivative of the reconstruction at $x = \frac{\Delta x}{2}$ is a linear combination of $q_{i-\frac12}$, $\bar q_i$ and $q_{i+\frac12}$, and by construction will be a derivative approximation exact for parabolae.

\begin{figure}
 \centering
 \includegraphics[width=0.8\textwidth]{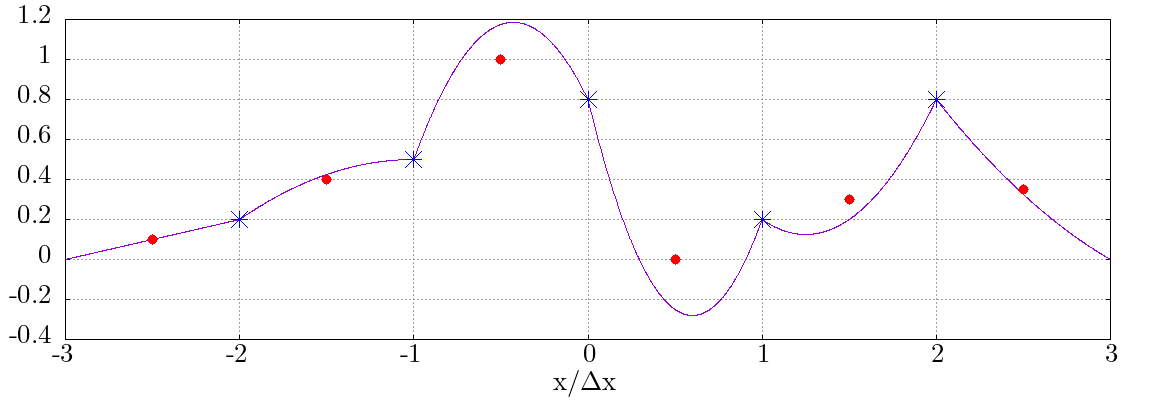}
 \caption{Example of a reconstruction/approximation of Active Flux for polynomial degree $K=2$. Point values are shown as crosses, averages as circles.}
 \label{fig:exampleglobal}
\end{figure}

By combining reconstructions in different cells one arrives at a global reconstruction
\begin{align}
 q_\text{recon} &\colon \mathbb R \to \mathbb R & q_\text{recon}(x) = q_{\text{recon},i}(x - x_i) \qquad \text{if }x\in [x_{i-\frac12}, x_{i+\frac12}]
\end{align}
Figure \ref{fig:exampleglobal} shows an example for such a global reconstruction.

In the Finite Element approach aimed at here this global reconstruction is also called the approximation or the numerical solution. It is piecewise polynomial and globally continuous, since
\begin{align}
 q_{\text{recon},i}\left( \frac{\Delta x}{2} \right ) = q_{\text{recon},i+1}\left(-\frac{\Delta x}{2} \right )
\end{align}

However, its derivative is not continuous at cell interfaces (see Figure \ref{fig:exampleglobal}). This is central to the way how upwinding is included in the semi-discrete Active Flux method. In Equation \eqref{eq:pointvaluesAF} one uses a left-/right-biased discrete derivative by differentiating either $q_{\text{recon},i}$ or $q_{\text{recon},i+1}$, i.e.
\begin{align}
 (D^+ q)_{i+\frac12} &:= \frac{\dd }{\dd x} q_{\text{recon},i}\left( \frac{\Delta x}{2}\right) & (D^- q)_{i+\frac12} &:= \frac{\dd }{\dd x} q_{\text{recon},i+1}\left(- \frac{\Delta x}{2}\right)
\end{align}
Through an increase in the size of the stencil one can achieve larger stability domains (as exemplified in \cite{abgrall22}), but this work focuses entirely on the approach of differentiating the reconstruction. The Active Flux method consists of the two equations \eqref{eq:averagesAF} and \eqref{eq:pointvaluesAF}, which can be solved with a Runge-Kutta method of correspondingly high order of accuracy. The method then is third-order accurate (see \cite{zeng14,barsukow19activeflux,kerkmann18} for an analyses of the accuracy of Active Flux).

It is possible to use other ways of finding discrete approximations to $\del_x f$, see \cite{duan24} for an example. It is also possible to generalize the Active Flux method to higher than third order. A particular case is the inclusion of higher moments, i.e. expressions such as $\int_{x_{i-\frac12}}^{x_{i+\frac12}} (x-x_i)^k q(t,x) \,\dd x$) as additional degrees of freedom (\cite{abgrall22}); see below for more details. In two spatial dimensions the procedure remains the same (with the derivatives perpendicular to edges being subject to upwinding and those parallel to the edges being unique and central). See \cite{barsukow24afeuler} and the corresponding Section \ref{sec:2d} below for further details.

The aim of the next section is to rephrase the approximation in the framework of Finite Elements and then to derive the Active Flux method from a Petrov-Galerkin variational framework.

\subsection{Degrees of freedom and basis functions}

Some concepts in Active Flux are related to Finite Elements in a very obvious way; it is useful to rephrase them first, starting with the Active Flux ``reconstruction''.
Active Flux uses a piecewise polynomial, globally continuous approximation in the sense of Finite Elements:

\begin{definition}[Active Flux Finite Element in 1D] \label{def:finele1d}
 Consider the interval $E = \left[- \frac{\Delta x}{2}, \frac{\Delta x}{2}\right]$ and denote by $V$ the space $P^K(E)$, $K \geq 2$. Denote by $V' = \mathrm{hom}(V, \mathbb R)$ the dual space of $V$ and by $(a_k)_{k=0, \ldots, K-2}$ a basis of $P^{K-2}$. Then, the \emph{Active Flux Finite Element} is the triple $(E, V, \Sigma)$ with $\Sigma \subset V'$ the following set of 
 degrees of freedom $\sigma_{\pm\frac12}, \sigma^{(k)} \in V'$ (point values and moments):
\begin{align}
 \sigma_{\pm\frac12}(v) &:= v\left(\pm \frac{\Delta x}{2}\right) &\sigma_k(v) &:= 
 \int_{E} \mombas_k(x) v(x) \, \dd x \qquad  \forall v \in V, k=0, \ldots, K-2
\end{align}

\end{definition}
For simplicity, we consider the polynomial degree $K$ fixed beforehand, and omit making explicit in the notation that $V$ and $\mathscr V$ (and other polynomial spaces introduced later) depend on $K$.

\begin{figure}[h]
 \centering
 \includegraphics[width=0.2\textwidth]{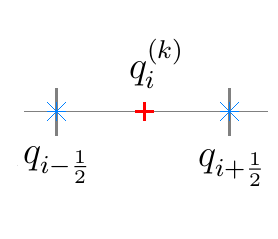}
 \caption{Degrees of freedom of Active Flux in 1-d.}
 \label{fig:dof1d}
\end{figure}

The classical choice (\cite{abgrall22}) is $\mombas_k(x) =\frac{(k+1)2^k}{\Delta x^{k+1}} x^k$, such that $\sigma_k(1) = 1 \, \forall k \in 2 \mathbb N_0$. 
The lowest-order moment for $\mombas_0(x) = 1$ is the cell average.
A choice of moments with respect to e.g. Legendre polynomials instead of the monomials $x^k$ is possible, but does not make any conceptual difference, and in particular does not seem to lead to any simplifications.
This Finite Element is essentially the hybrid Finite Element given in Prop. 6.10 (Section 6.3.3, p. 66) in \cite{guermondfinele1}. An illustration of the degrees of freedom is given in Figure \ref{fig:dof1d}.

\begin{figure}
 \centering
 \includegraphics[width=0.28\textwidth]{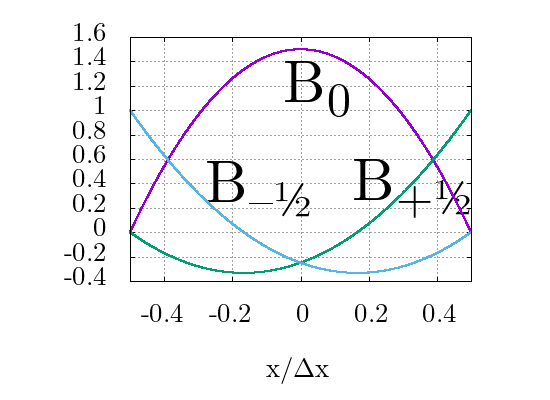} \hfill \includegraphics[width=0.64\textwidth]{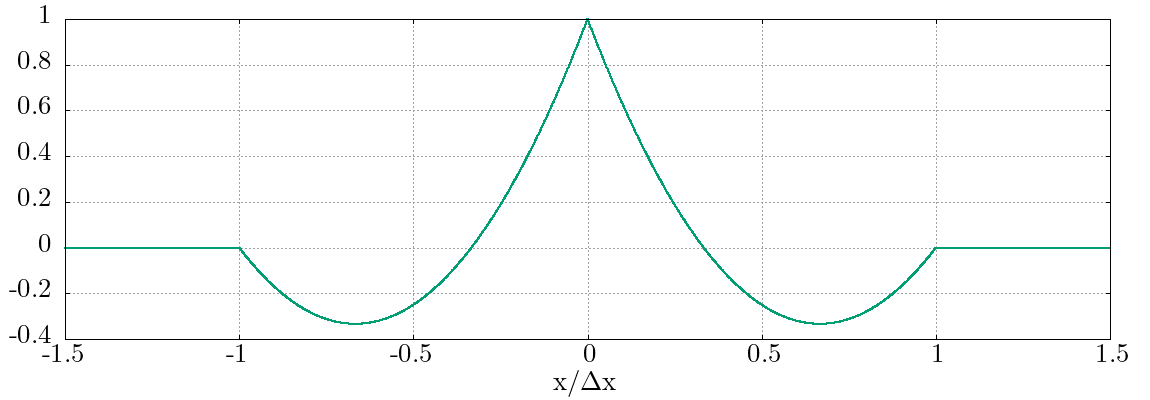}
 \caption{The basis functions for $K=2$. \emph{Left}: The (local) basis functions $B_0, B_{\pm\frac12} \in V$. \emph{Right}: The global basis function $\phi_{i+\frac12} \in \mathscr V$.}
 \label{fig:basis}
\end{figure}

The space $V$ can be endowed with a basis $\{ B_{\pm\frac12} \} \cup \{ B_k \}_{k = 0, \dots, K-2}$ dual to $\Sigma$, i.e.
\begin{align}
 \sigma_r(B_s) = \delta_{rs} \qquad \forall r,s \in \left \{ \pm\frac12, 0, 1, \dots, K-2 \right \} \label{orthogonality1d}
\end{align}
$B_{\pm\frac12}$ are the basis functions associated to the left/right point value, and $B_k$ is the basis function associated to the $k$-th moment (see Figure \ref{fig:basis} for an example).

\begin{example}
 For $K=2$, one finds, with $\xi := \frac{x}{\Delta x}$,
 \begin{align}
 B_0(x) &= \frac{3}{2} \left( 1 - 4 \xi^2 \right ) &
 B_{+\frac12}(x) &= \frac14 \left( 2\xi+1\right ) \left( 6  \xi- 1\right ) 
 \end{align}
 \begin{align}
 B_{-\frac12}(x) &= \frac14 \left( 2 \xi-1\right ) \left( 6  \xi + 1\right )
\end{align}
\end{example}

Next, a switch from the reference element to the global view is needed. Here, it is important to take into account that the point values at cell interfaces are identified, as usual in continuous Finite Elements. The corresponding global basis function thus has support in two elements, while the global basis functions associated to the moments continue to be supported in just one element.

\begin{definition}
 Consider the space
 \begin{align}
 \mathscr V := \{ v : v \in C^0 \cap L^2, v\vert _{C_i} \in V \, \forall i \in \mathbb Z\} 
\end{align}
Then, denote by $(\phi_i^{(k)})_{i\in\mathbb Z, k = 0, \dots, K-2} \cup (\phi_{i+\frac12})_{i\in \mathbb Z}$ the \emph{global basis} of $\mathscr V$ induced by the basis $\{ B_{\pm\frac12} \} \cup \{ B_k \}_{k = 0, \dots, K-2}$ of $V$:
\begin{align}
 \phi_i^{(k)}(x) &:= B_k(x - x_i) \chi_{[x_{i-\frac12},x_{i+\frac12}]}\\
 \phi_{i+\frac12}(x) &:= B_{+\frac12}(x - x_i) \chi_{[x_{i-\frac12},x_{i+\frac12})} + B_{-\frac12}(x - x_{i+1}) \chi_{[x_{i+\frac12},x_{i+\frac32}]}
\end{align}
\end{definition}
See Figure \ref{fig:basis} for an illustration. In presence of boundaries the infinite grid needs to be replaced by a finite one. At the domain boundaries, it is customary for Active Flux to use ghost cells. A special case is a Dirichlet boundary coinciding with cell interfaces, in which case the point values at cell interfaces can be imposed directly and are not solved for (see e.g. \cite{barsukow23afim}). A more detailed discussion of boundary conditions is beyond the scope of this work, however.

This allows to write a numerical solution $q_h(t, \cdot) \in \mathscr V$ as
\begin{align}
 q_h(t, x) = \sum_{i \in \mathbb Z} \left( \phi_{i+\frac12}(x) q_{i+\frac12}(t) + \sum_{k = 0}^{K-2} \phi^{(k)}_i(x) q^{(k)}_i(t)   \right ) \label{eq:numapprox}
\end{align}

The orthogonality \eqref{orthogonality1d} is now replaced by the following relations, valid for all $i,j \in \mathbb Z$, $k,\ell=0, \ldots, K-2$:
\begin{align}
 \phi_{i+\frac12}(x_{j+\frac12}) &= \delta_{ij} & \phi^{(k)}_i(x_{j+\frac12}) = 0  \label{eq:generalbasiscondition1}\\
 \int_{x_{j-\frac12}}^{x_{j+\frac12}} \mombas_k(x-x_i) \phi_{j+\frac12}(x) \,\dd x &= 0 & 
 \int_{x_{j-\frac12}}^{x_{j+\frac12}} \mombas_\ell(x-x_i) \phi^{(k)}_i(x) \,\dd x &= \delta_{k\ell} \delta_{ij} \label{eq:generalbasiscondition2}
\end{align}
where, again, the classical choice adopted here is $\mombas_k$ proportional to $x^k$ as in Definition \ref{def:finele1d}.

\subsection{Biorthogonal test functions and the Petrov-Galerkin method} \label{sec:biorthogonal}

Bi\-or\-thogonal polynomials are a generalization of families of orthogonal polynomials, and usually the degree of each family is imposed to be monotonically increasing. \cite{konhauser65} gives conditions for the existence of such biorthogonal families.
The degree of nodal basis functions is usually the same, and thus biorthogonality in the context of Finite Elements leads to a slightly different problem. Dual or biorthogonal test functions appear e.g. in \cite{tchamitchian87,scott90,goodman00,wohlmuth00,oswald01,lamichhane07,dickopf14,haubold24} in various contexts, with wavelet transforms being among the earliest. They do not seem to be used very often for classical Finite Element methods; this work aims at demonstrating that semi-discrete Active Flux arises naturally as a Petrov-Galerkin method with biorthogonal test functions.

We seek a set of test functions $\big( \psi^{(k)}_i  \big)_{i\in\mathbb Z, k= 0, \ldots, K-2}  \cup \big( \psi_{i+\frac12} \big)_{i\in\mathbb Z}  $ with the following\footnote{With a global counting of degrees of freedom one would simply have $(\psi_{\alpha}, \phi_{\beta})_{L^2} = \delta_{\alpha\beta}$. Here, due to the different nature of point values and averages, it is preferable to distinguish the statements.} orthogonality properties:
\begin{align}
 (\psi_{j+\frac12}, \phi_{i+\frac12})_{L^2} &= \delta_{ij} & (\psi_{j+\frac12}, \phi^{(k)}_i)_{L^2} &= 0 \label{eq:biorthopoint}\\
 (\psi^{(\ell)}_j, \phi_{i+\frac12})_{L^2} &= 0 & (\psi^{(\ell)}_j, \phi^{(k)}_i)_{L^2} &= \delta_{k\ell} \delta_{ij} \label{eq:biorthomoment}
\end{align}

Define $\mathscr W := \mathrm{span}_{i \in \mathbb Z, \,k= 0, \ldots, K-2}(\psi^{(k)}_i, \psi_{i+\frac12} )$. As a short-hand, they will be called ``biorthogonal test functions'' even though it would be more appropriate to speak of a biorthogonal pair of piecewise-polynomial families (formed by the basis functions $\phi$ and the test functions $\psi$).

\begin{definition}
 Given the numerical approximation \eqref{eq:numapprox} and having fixed the polynomial degree $K$, the \emph{Active Flux Petrov-Galerkin method} for the conservation law \eqref{eq:conslaw1d}, i.e. $\del_t q + \del_x f(q) = 0$ is given by
 \begin{align}
  (\psi, \del_t q_h)_{L^2} + (\psi, \del_x f(q_h))_{L^2} &= 0 \qquad \forall \psi \in \mathscr W
 \end{align}
\end{definition}

In particular the method amounts to the two sets of equations
\begin{align}
  \frac{\dd }{\dd t} q_i^{(k)}  + (\psi^{(k)}_i, \del_x f(q_h))_{L^2} &= 0 \qquad \forall i \label{eq:pgmom1d}\\
  \frac{\dd}{\dd t} q_{i+\frac12} + (\psi_{i+\frac12}, \del_x f(q_h))_{L^2} &= 0 \qquad \forall i \label{eq:pgpoint1d}
 \end{align}
where the diagonal mass matrix is a consequence of the biorthogonality \eqref{eq:biorthopoint}--\eqref{eq:biorthomoment}.

\begin{theorem} \label{thm:testmoment1d}
 The test functions associated to the moments are $\psi^{(k)}_i =  \mombas_k(x - x_i) \id_{[x_{i-\frac12}, x_{i+\frac12}]}$. 
\end{theorem}
\begin{proof}
 Indeed, for all $j \in \mathbb Z$,
 \begin{align}
  (\psi^{(\ell)}_i, \phi_{j+\frac12})_{L^2} &= \int_{\mathbb R} 
  \mombas_\ell(x - x_i) \id_{[x_{i-\frac12}, x_{i+\frac12}]} \phi_{j+\frac12}(x) \,\dd x = 0 \label{eq:deftestfct1d}\\
  (\psi^{(\ell)}_i, \phi^{(k)}_j)_{L^2} &= \int_{\mathbb R} 
  \mombas_\ell(x - x_i) \id_{[x_{i-\frac12}, x_{i+\frac12}]} \phi^{(k)}_j(x) \,\dd x =\delta_{k\ell} \delta_{ij}
 \end{align}
 by definition.
\end{proof}

Theorem \ref{thm:testmoment1d} implies that the evolution equations \eqref{eq:pgmom1d} for the moments are just those of Active Flux (as introduced in \cite{abgrall22,barsukow24afeuler}):
\begin{align}
 \frac{\dd }{\dd t} q_i^{(k)} + 
 \int_{x_{i-\frac12}}^{x_{i+\frac12}}\mombas_k(x-x_i) \del_x f(q_h) \, \dd x = 0
\end{align}
To proceed, partial integration is usually applied to the integral now, see \cite{abgrall22} for further details. In particular, for the cell averages one obtains
\begin{align}
 \frac{\dd }{\dd t} q_i^{(0)} + \frac{f(q_{i+\frac12}) - f(q_{i-\frac12})}{\Delta x} = 0
\end{align}

\begin{theorem}
 There exists a 1-parameter family of functions $\psi_{i+\frac12}$ supported in $[x_{i-\frac12}, x_{i+\frac32}]$ with restrictions onto $(x_{i-\frac12},x_{i+\frac12})$ and $(x_{i+\frac12},x_{i+\frac32})$ in $P^K$ which solve \eqref{eq:biorthopoint}.
\end{theorem}
\begin{proof}
Assume that $\overline{\mathrm{supp}\,\psi_{i+\frac12}} = \overline{\mathrm{supp}\,\phi_{i+\frac12}} = [x_{i-\frac12}, x_{i+\frac32}]$ and denote by
\begin{align}
 \phi^{\text{L}}_{i+\frac12}(x) := \phi_{i+\frac12}(x) \Big \vert _{[x_{i-\frac12}, x_{i+\frac12}]} = B_{+\frac12}(x-x_i)\\
 \phi^{\text{R}}_{i+\frac12}(x) := \phi_{i+\frac12}(x) \Big \vert _{[x_{i+\frac12}, x_{i+\frac32}]} = B_{-\frac12}(x - x_{i+1})
\end{align}
and similarly 
\begin{align}
 \psi^{\text{L}}_{i+\frac12}(x) := \psi_{i+\frac12}(x) \Big \vert _{[x_{i-\frac12}, x_{i+\frac12}]} =: A_{+\frac12}(x-x_i)\\
 \psi^{\text{R}}_{i+\frac12}(x) := \psi_{i+\frac12}(x) \Big \vert _{[x_{i+\frac12}, x_{i+\frac32}]} =: A_{-\frac12}(x - x_{i+1})
\end{align}
$A_{\pm\frac12} \colon [-\frac{\Delta x}{2} , \frac{\Delta x}{2}] \to \mathbb R$ is thus the local form of the test function. 
The defining relations \eqref{eq:biorthopoint} of $\psi_{j+\frac12}$ therefore are
\begin{align}
 1&= \int_{x_{i-\frac12}}^{x_{i+\frac12}} \psi_{i+\frac12}^{\text{L}} \phi_{i+\frac12}^{\text{L}} \, \dd x +\int_{x_{i+\frac12}}^{x_{i+\frac32}} \psi_{i+\frac12}^{\text{R}} \phi_{i+\frac12}^{\text{R}} \, \dd x \\
 &\qquad\qquad\qquad= \int_{-\frac{\Delta x}{2}}^{\frac{\Delta x}{2}} A_{+\frac12} B_{+\frac12} \, \dd x + \int_{-\frac{\Delta x}{2}}^{\frac{\Delta x}{2}} A_{-\frac12} B_{-\frac12} \, \dd x \label{eq:testfct1}\\
 0 &= \int_{x_{i-\frac12}}^{x_{i+\frac12}} \psi_{i-\frac12}^{\text{R}} \phi_{i+\frac12}^{\text{L}} \, \dd x = \int_{-\frac{\Delta x}{2}}^{\frac{\Delta x}{2}} A_{-\frac12} B_{+\frac12} \, \dd x \label{eq:testfct2}\\
 0 &= \int_{x_{i-\frac12}}^{x_{i+\frac12}} \psi_{i+\frac12}^{\text{L}} \phi_{i-\frac12}^{\text{R}} \, \dd x = \int_{-\frac{\Delta x}{2}}^{\frac{\Delta x}{2}} A_{+\frac12} B_{-\frac12} \, \dd x\\
 0 &= \int_{x_{i-\frac12}}^{x_{i+\frac12}} \psi_{i+\frac12}^{\text{L}} \phi^{(k)}_{i} \, \dd x = \int_{-\frac{\Delta x}{2}}^{\frac{\Delta x}{2}} A_{+\frac12} B_k \, \dd x \quad \forall k = 0, \ldots, K-2\\
 0 &= \int_{x_{i-\frac12}}^{x_{i+\frac12}} \psi_{i-\frac12}^{\text{R}} \phi^{(k)}_{i} \, \dd x = \int_{-\frac{\Delta x}{2}}^{\frac{\Delta x}{2}} A_{-\frac12} B_k \, \dd x\quad \forall k = 0, \ldots, K-2 \label{eq:testfctlast}
\end{align}
These are $2(K-1)+3 = 2K+1$ equations for the two functions $A_{\pm\frac12}$. If these shall be taken from $P^{K}([-\frac{\Delta x}{2}, \frac{\Delta x}{2}])$, then there are a total of $2(K+1)$ parameters to be determined, i.e. one parameter is going to remain free. It is chosen here as follows: a parameter $\alpha \in \mathbb R$ is introduced such that \eqref{eq:testfct1} is replaced by
\begin{align}
\frac12 + \frac{\alpha}{2} &= \int_{-\frac{\Delta x}{2}}^{\frac{\Delta x}{2}} A_{+\frac12} B_{+\frac12} \, \dd x \label{eq:centraltestfct1}\\
\frac12 - \frac{\alpha}{2} &= \int_{-\frac{\Delta x}{2}}^{\frac{\Delta x}{2}} A_{-\frac12} B_{-\frac12} \, \dd x \label{eq:centraltestfct2}
\end{align}

The linear systems for $A_{\pm\frac12}$ have a unique solution.  
This is established as follows: Focus on those Equations \eqref{eq:testfct2}--\eqref{eq:testfctlast}, \eqref{eq:centraltestfct1}--\eqref{eq:centraltestfct2} which concern $A_{+\frac12}$, for definiteness. Since $\{ B_{\pm\frac12}, B_0, \dots, B_{K-2} \}$ form a basis of $P^K([-\frac{\Delta x}{2}, \frac{\Delta x}{2}])$, all Legendre polynomials $L_s$ of degree at most $K$, which span $P^K([-\frac{\Delta x}{2}, \frac{\Delta x}{2}])$, can be expressed as linear combinations of $\{ B_{\pm\frac12}, B_0, \dots, B_{K-2} \}$. Thus, the linear system becomes
\begin{align}
 \int_{-\frac{\Delta x}{2}}^{\frac{\Delta x}{2}} A_{+\frac12} L_s \, \dd x = \textrm{RHS} \qquad s = 0, \dots, K
\end{align}
Now, one can also consider $A_{+\frac12} \in P^K$ in the Legendre basis as $A_{+\frac12} = \sum_{r = 0}^{K} a_{+\frac12,r} L_r$, as well. But then the linear system is simply
\begin{align}
 \sum_{r = 0}^{K} a_{+\frac12,r}\int_{-\frac{\Delta x}{2}}^{\frac{\Delta x}{2}}  L_r L_s \, \dd x = \textrm{RHS} \qquad s = 0, \dots, K
\end{align}
Thus, due to $L^2$-orthogonality of Legendre polynomials, in this basis, the linear systems are diagonal and the entries are non-zero.
\end{proof}

\begin{figure}
 \centering
 \includegraphics[width=0.28\textwidth]{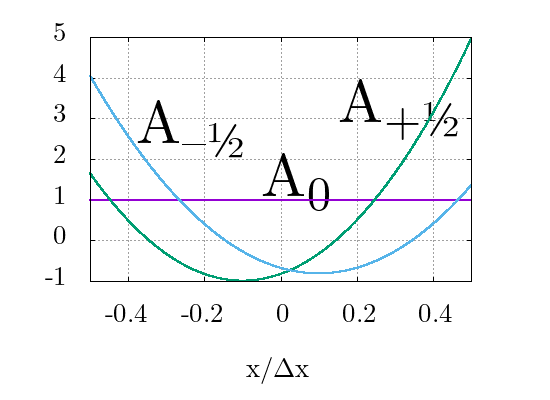} \hfill \includegraphics[width=0.64\textwidth]{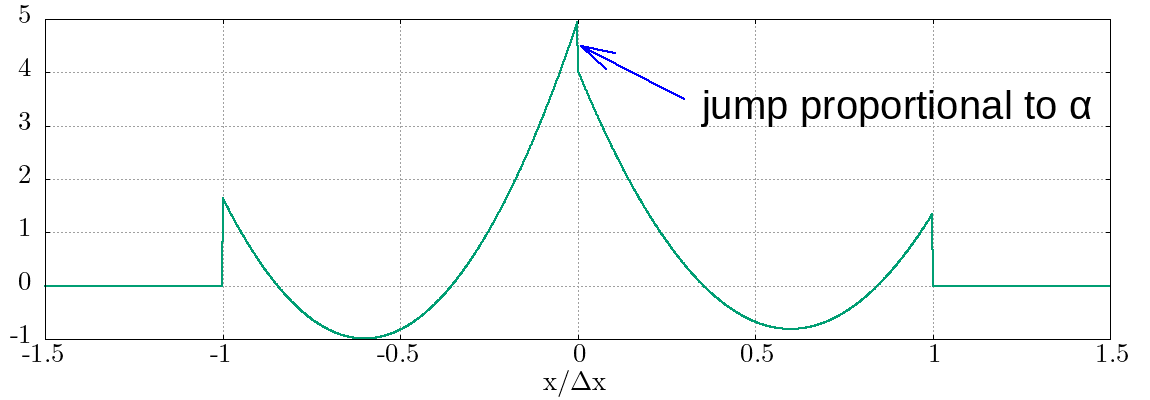} \\
 \includegraphics[width=0.44\textwidth]{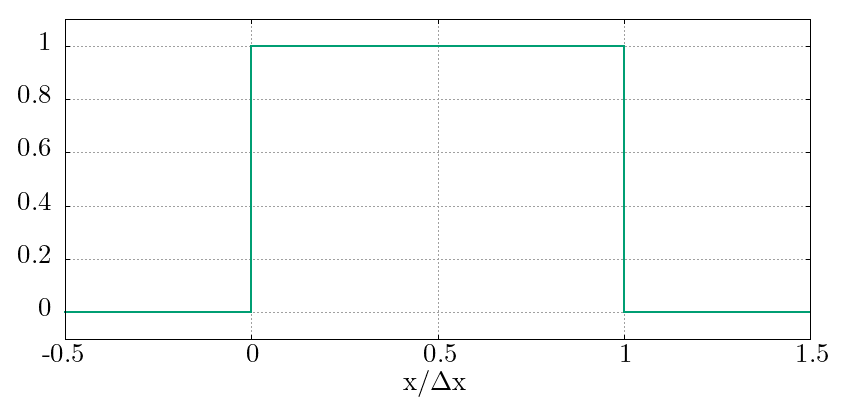}
 \caption{The test functions for $K=2$. \emph{Left}: The (local) test functions $A_0, A_{\pm\frac12}$. \emph{Right}: The global test function $\psi_{i+\frac12} $. \emph{Bottom}: The global test function $\psi_{i}^{(0)}$. Here, $\alpha = 0.1$.} 
 \label{fig:test}
\end{figure}

\begin{example}
 For $K=2$ one finds, with $\alpha \in \mathbb R$ and $\xi := x/\Delta x$,
 \begin{align}
  A_{\pm\frac12}(x) &= \frac34 (1 \pm \alpha) \left( -1 \pm 4 \xi + 20 \xi^2 \right)
 \end{align}
  and in particular, if $\alpha = 0$
 \begin{align}
  A_{\pm\frac12}(x) &= -\frac34  \pm 3 \xi + 15 \xi^2 \label{eq:1dbiorthtestsymm}
 \end{align}
 
 Centered at $x_{i+\frac12}$, one finds
 \begin{align}
  \psi_{i+\frac12}(x) &= \frac{3}{2} \left(1 - \alpha_{i+\frac12} \, \mathrm{sgn}(x - x_{i+\frac12})\right) \left(  3 - 12 \frac{\vert x - x_{i+\frac12}\vert }{\Delta x} + 10 \frac{(x - x_{i+\frac12})^2}{\Delta x^2}   \right)
 \end{align}
 where the free parameter $\alpha_{i+\frac12}$ can be chosen differently at each interface.
 See also Figure \ref{fig:test}.
\end{example}

In general, the test function $\psi_{i+\frac12}$ is discontinuous at $x_{i-\frac12}$, $x_{i+\frac12}$, and $x_{i+\frac32}$, with the jump at $x_{i+\frac12}$ proportional to $\alpha_{i+\frac12}$ and those at the endpoints of the support proportional to $1 \pm \alpha_{i+\frac12}$. The discontinuity at $x_{i+\frac12}$ allows to control the amount of upwinding present in the discrete derivative $(\psi_{i+\frac12},\del_x q(t,x))_{L^2}$, as demonstrated by the following

\begin{theorem} \label{thm:derivatives}
 Consider the Petrov-Galerkin method \eqref{eq:pgpoint1d} for linear $f$. Then, $(Dq)_{i+\frac12} := (\psi_{i+\frac12},\del_x q_h(t,x))_{L^2}$ is the $\alpha$-upwinded average
 \begin{align}
  (Dq_h)_{i+\frac12} = \frac{d_{i+\frac12}^+  + d_{i+\frac12}^- }{2} - \frac{\alpha_{i+\frac12}}{2} (d_{i+\frac12}^- - d_{i+\frac12}^+)
 \end{align}
 of the derivatives 
 \begin{align}
 d_{i+\frac12}^+ &:= \lim_{\epsilon \searrow 0} (\del_x q_h)(t,x_{i+\frac12} - \epsilon ) &
 d_{i+\frac12}^- &:= \lim_{\epsilon \searrow 0} (\del_x q_h)(t,x_{i+\frac12} + \epsilon )  \label{eq:derivativedef}
 \end{align}
 \end{theorem}
 \begin{proof}
 The idea is to expand the derivative $\del_x q_h$ in each cell in the Finite Element basis.
 Since $P^K([x_{i-\frac12},x_{i+\frac12}]) = \mathrm{span}\left(\phi_{i-\frac12}^{\text{R}}, \phi^{(k)}_{i}, \phi_{i+\frac12}^{\text{L}}\right)$, and since $\del_x q_h \Big\vert _{(x_{i-\frac12},x_{i+\frac12})} \in P^K$, there exist expansion coefficients $d_{i-\frac12}^-, d^{(0)}_i, \dots, d^{(K-2)}_i, d_{i+\frac12}^+$ such that
 \begin{align}
  \del_x q_h \Big\vert _{(x_{i-\frac12},x_{i+\frac12})} &= \del_x \phi_{i-\frac12}^{\text{R}} q_{i-\frac12} + \sum_{k = 0}^K \del_x \phi_i^{(k)}(x) q_i^{(k)} +\del_x \phi_{i+\frac12}^{\text{L}} q_{i+\frac12} \\
  &= d_{i-\frac12}^- \phi_{i-\frac12}^{\text{R}} + \sum_{k = 0}^K d^{(k)}_i \phi^{(k)}_{i} + d_{i+\frac12}^+ \phi_{i+\frac12}^{\text{L}}
 \end{align}
 where explicit knowledge of $d^{(k)}_i$ is unnecessary for what follows. It is, however, clear from the definition of the basis functions that indeed, $d_{i-\frac12}^-, d_{i+\frac12}^+$ have the meaning \eqref{eq:derivativedef}.
  Then,
 \begin{align}
 \int_{\mathbb R} \psi_{i+\frac12} &\del_x q_h(t,x) \, \dd x = \int_{x_{i-\frac12}}^{x_{i+\frac12}} \psi_{i+\frac12}^{\text{L}} \left(  d_{i-\frac12}^- \phi_{i-\frac12}^{\text{R}} + \sum_{k = 0}^{K-2} d^{(k)}_i \phi^{(k)}_{i} + d_{i+\frac12}^+ \phi_{i+\frac12}^{\text{L}} \right ) \dd x\\ 
 &+ \int_{x_{i+\frac12}}^{x_{i+\frac32}} \psi_{i+\frac12}^{\text{R}} \left(  d_{i+\frac12}^- \phi_{i+\frac12}^{\text{R}} + \sum_{k = 0}^{K-2} d^{(k)}_{i+1} \phi^{(k)}_{i+1} + d_{i+\frac32}^+ \phi_{i+\frac32}^{\text{L}} \right ) \dd x \\
 &\!\!\!\!\!\!\!\overset{\text{\eqref{eq:testfct1}--\eqref{eq:testfctlast}}}{=}  d_{i+\frac12}^+ \int_{x_{i-\frac12}}^{x_{i+\frac12}} \psi_{i+\frac12}^{\text{L}}  \phi_{i+\frac12}^{\text{L}}  \dd x +  d_{i+\frac12}^- \int_{x_{i+\frac12}}^{x_{i+\frac32}} \psi_{i+\frac12}^{\text{R}}  \phi_{i+\frac12}^{\text{R}}  \dd x \\
 &\!\!\!\!\!\!\!\overset{\text{\eqref{eq:centraltestfct1}--\eqref{eq:centraltestfct2}}}{=} \frac12 d_{i+\frac12}^+ \left( 1 + \alpha \right) + \frac12 d_{i+\frac12}^- \left( 1 - \alpha   \right ) \\
 &= \frac{d_{i+\frac12}^+  + d_{i+\frac12}^- }{2} - \frac12 \alpha (d_{i+\frac12}^- - d_{i+\frac12}^+) =
 \begin{cases}
  d_{i+\frac12}^+ & \alpha = 1 \\
  d_{i+\frac12}^- & \alpha = -1
 \end{cases}
 \end{align}
\end{proof}

\begin{figure}
 \centering
 \includegraphics[width=0.49\textwidth]{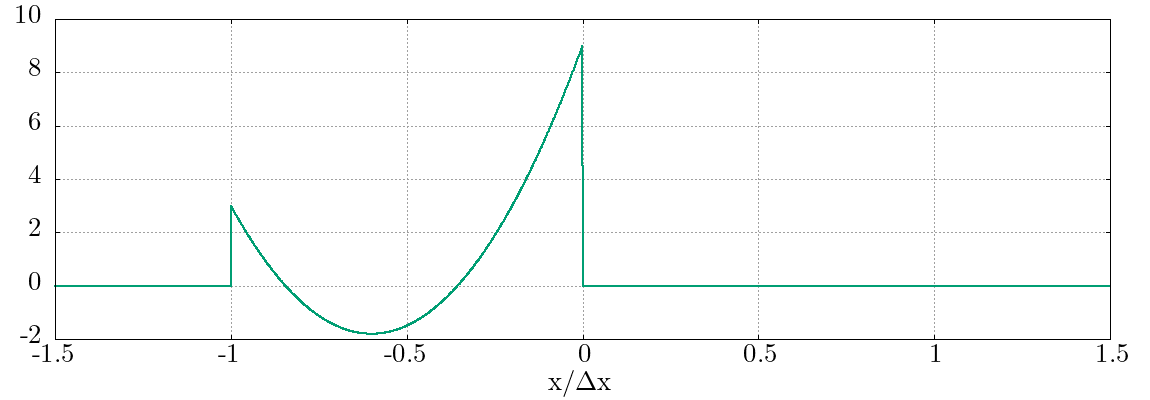} \hfill \includegraphics[width=0.49\textwidth]{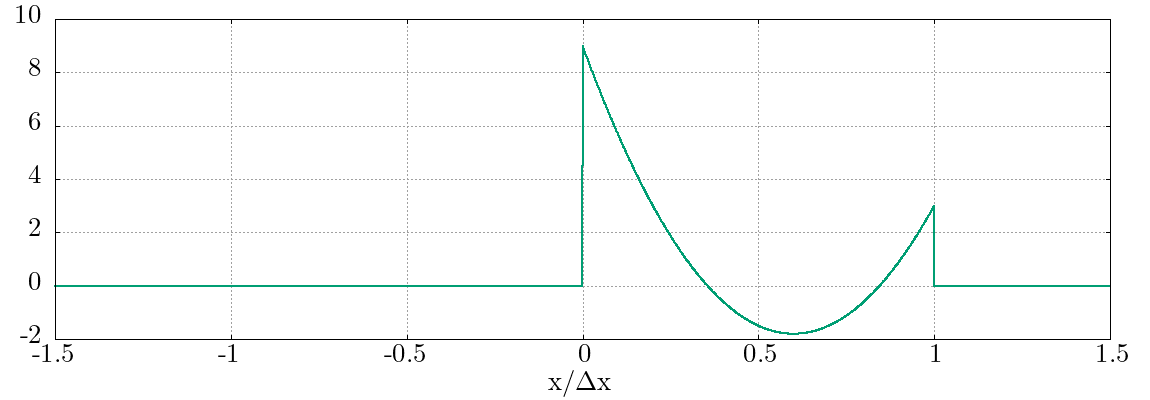} \\
 \caption{The global test function $\psi_{i+\frac12}$ for $K=2$ at $x_{i+\frac12} = 0$. \emph{Left}: $\alpha_{i+\frac12} = 1$. \emph{Right}: $\alpha_{i+\frac12} = -1$. Here, $\Delta x = 1$.} 
 \label{fig:testOne}
\end{figure}

The Theorem thus shows that the free parameter $\alpha$ is just the usual parameter governing the amount of numerical diffusion in the upwind numerical flux, if $f$ is linear. For linear problems, the upwind direction is predetermined, while for nonlinear ones it can change in space and time, such that $\alpha$ needs to be chosen adaptively.
One thus in general would want to choose $\alpha_{i+\frac12}$ as $\mathrm{sgn} f'(q_{i+\frac12})$ by analogy with upwinding based on Riemann solvers and linearization. Independently of whether $f$ is linear or nonlinear, $\alpha$ has an influence on the support of the test function: the test function is fully supported in the upwind or the downwind cell for values $\alpha = \pm1$ (see Figure \ref{fig:testOne}), which is what upwinding is about. The jumps of the derivatives of the approximation vanish if the solution is globally polynomial of degree at most $K$, i.e. $\alpha$ does not change the order of accuracy of the method.

The Streamline-Upwind (SUPG) Petrov-Galerkin method is another example of upwinding introduced through a discontinuity in the test function (SUPG involves derivatives of the basis functions in its test functions).

The above method is precisely the Active Flux method from \cite{abgrall22} in the cases of full upwinding, i.e. $\alpha = \pm1$, 
 \begin{align}
  \Delta x D &= 2 q_{i-\frac12} -6 q_{i} +4 q_{i+\frac12} \\
  \Delta x D^* &= -4 q_{i+\frac12} +6 q_{i+1} -2 q_{i+\frac32} 
 \end{align}
while the average (i.e. $\alpha = 0$) is the central Active Flux 
 \begin{align}
  \left(\frac12 (D+D^*) q\right)_{i+\frac12} &= \frac{q_{i-\frac12} - 3 q_i + 3 q_{i+1} - q_{i+\frac32}}{\Delta x} \label{eq:derivativecentral}
 \end{align}
considered e.g. in \cite{barsukow24affourier} and whose energy stability is analyzed in \cite{barsukow2025sbp}.

For nonlinear problems, one can either evaluate $(\psi_{i+\frac12}, \del_x f(q_h))_{L^2}$ exactly, or project $f(q_h)$ onto $\mathscr V$ first. The latter approach leads to the flux-vector splitting introduced in \cite{duan24}. An example of the former is the following.

\begin{example}
 For Burgers' equation with $f(q) = \frac{q^2}{2}$ one finds
 \begin{align}
  (\psi_{i+\frac12}, \del_x f(q_h))_{L^2} &= \frac{1+\alpha_{i+\frac12}}{2} \cdot \frac{-9 (q_{i-\frac12} -2 q_i)^2 + 2 (q_{i-\frac12}-12 q_i)q_{i+\frac12} + 31 q_{i+\frac12}^2}{10 \Delta x} \label{eq:derivburgers} \\&\!\!\!\!\!\!\!\!\!\!\!\!\!\!+\nonumber
  \frac{1-\alpha_{i+\frac12}}{2} \cdot \frac{9 (q_{i+\frac32} -2 q_{i+1})^2 - 2 (q_{i+\frac32}-12 q_{i+1})q_{i+\frac12} - 31 q_{i+\frac12}^2}{10 \Delta x}
 \end{align}
 Experimentally, this discretization does not perform differently from the choices in \cite{abgrall22,duan24} based on a linearization or flux-vector splitting, see Section \ref{sec:numericalburgers}.
\end{example}

\section{Two spatial dimensions} \label{sec:2d}

\subsection{Degrees of freedom and basis functions}

\begin{figure}[h]
 \centering
 \includegraphics[width=0.35\textwidth]{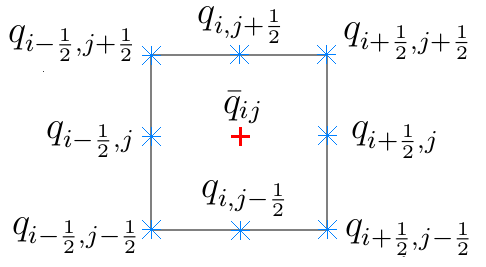}
 \caption{Degrees of freedom of Active Flux in 2-d for $K=2$.}
 \label{fig:dof2d}
\end{figure}

For the sake of simplicity we restrict ourselves here to the case $K=2$, i.e. the biparabolic one. The degrees of freedom are the cell average $\bar q_{ij}$, as well as the point values at edge midpoints $q_{i+\frac12,j}$, $q_{i,j+\frac12}$ and at nodes $q_{i+\frac12,j+\frac12}$ (see Figure \ref{fig:dof2d}). Observe that this choice does not give rise to a tensor-product structure of the corresponding basis functions.

Clearly, it is possible to define basis functions both locally and globally in complete analogy to the one-dimensional case. In every cell it is possible to construct a biparabolic polynomial, traditionally referred to as reconstruction. The expressions for the local basis functions are 
\begin{align}
B_{+\frac12,+\frac12}(x,y) &= \frac{1}{16} (2 \xi+1) (2 \eta+1) (-1+2 \eta+2 \xi +12\xi \eta) \label{eq:basis2d1}\\
B_{-\frac12,+\frac12}(x,y) &=  \frac{1}{16} (2 \xi-1) (2 \eta+1) (1-2 \eta+2 \xi +12\xi \eta)\\
B_{-\frac12,-\frac12}(x,y) &=  \frac{1}{16} (2 \xi-1) (2 \eta-1) (-1-2 \eta-2 \xi +12\xi \eta) \\
B_{+\frac12,-\frac12}(x,y) &=  \frac{1}{16} (2 \xi+1) (2 \eta-1) (1+2 \eta-2 \xi +12 \xi\eta)\\
B_{0,j+\frac12}(x,y) &=  -\frac{1}{4} (4 \xi^2-1) (2 \eta+1) (6 \eta-1) \\
B_{0,j-\frac12}(x,y) &=  -\frac{1}{4} (4 \xi^2-1) (2 \eta-1) (6 \eta+1) \\
B_{-\frac12,0}(x,y) &=  -\frac{1}{4} (2 \xi-1) (6 \xi+1) (4 \eta^2-1) \\ 
B_{+\frac12,0}(x,y) &=  -\frac{1}{4} (2 \xi+1) (6 \xi-1) (4 \eta^2-1) \\
 B_0(x,y) &= \frac{9}{4} (4 \xi^2 -1) (4 \eta^2 -1) \label{eq:basis2dlast}
\end{align}
with $\xi = x/\Delta x$ and $\eta = y/\Delta y$.
We denote by $B_0$ the basis function associated to the average, by $B_{\pm\frac12,\pm\frac12}$ those associated to the nodes and by $B_{\pm\frac12,0}, B_{0,\pm\frac12}$ those of the edge midpoints:
\begin{align}
 B_0, B_{\pm\frac12,\pm\frac12}, B_{\pm\frac12,0}, B_{0,\pm\frac12} \colon \left[ -\frac{\Delta x}{2}, \frac{\Delta x}{2} \right ] \times \left[ -\frac{\Delta y}{2}, \frac{\Delta y}{2} \right ] \to \mathbb R
\end{align}
Since the edges are shared by two adjacent cells, and nodes by four, the numerical solution can be written as 
\begin{align}
 q_h(t,x,y) &= \sum_{(i,j) \in \mathbb Z^2} \Big ( \phi_{i+\frac12,j+\frac12}(x,y) q_{i+\frac12,j+\frac12}(t) + \phi_{i+\frac12,j}(x,y) q_{i+\frac12,j}(t)  \\ \nonumber & \phantom{mmmmmmmmmmmmm}+ \phi_{i,j+\frac12}(x,y) q_{i,j+\frac12}(t) + \phi_{ij}(x,y) \bar q_{ij}(t)   \Big) \label{eq:numapprox2d}
\end{align}
having introduced global basis functions
\begin{align}
 \phi_{i+\frac12,j+\frac12}(x,y) &= B_{+\frac12,+\frac12}(x - x_i,y-y_j) \id_{C_{ij}} + B_{+\frac12,-\frac12}(x - x_i,y-y_{j+1}) \id_{C_{i,j+1}}
 \\ \nonumber &\!\!\!\!\!\!\!\!\!\!\!\!\!\!\!\!\!\!\!\!\!\!+ B_{-\frac12,+\frac12}(x - x_{i+1},y-y_j) \id_{C_{i+1,j}} + B_{-\frac12,-\frac12}(x - x_{i+1},y-y_{j+1}) \id_{C_{i+1,j+1}} \\
 \phi_{i+\frac12,0}(x,y) &= B_{+\frac12,0}(x - x_i,y-y_j) \id_{C_{ij}} + B_{-\frac12,0}(x - x_{i+1},y-y_j) \id_{C_{i+1,j}} \\ 
 \phi_{i,j+\frac12}(x,y) &= B_{0,+\frac12}(x - x_i,y-y_j) \id_{C_{ij}} + B_{0,-\frac12}(x - x_i,y-y_{j+1}) \id_{C_{i,j+1}}
\end{align}

Consider the conservation law \eqref{eq:conslaw2d}.
The update equation for the cell average is obtained by integrating over the cell, using Gauss' law and replacing the integrals along the edges by quadratures that employ the point values available:
\begin{align}
  \frac{\dd}{\dd t} \bar q_{ij} &+ \frac{1}{\Delta y} \int_{y_{j-\frac12}}^{y_{j+\frac12}} \frac{f(q_h(t, x_{i+\frac12}, y))  - f(q_h(t, x_{i-\frac12}, y))}{\Delta x} \dd y \label{eq:averagesAF2d}\\ \nonumber
  &+ \frac{1}{\Delta x} \int_{x_{i-\frac12}}^{x_{i+\frac12}} \frac{g(q_h(t, x, y_{j+\frac12}) ) - g(q_h(t, x, y_{j-\frac12}))}{\Delta y} \dd x = 0  
\end{align}
The natural quadrature here is Simpson's rule.

For the update of the point values it has been suggested in \cite{barsukow24afeuler} to use again the quasi-linear form of the conservation law:
\begin{align}
 \frac{\dd}{\dd t} q_p + J (D_x q)_p + K (D_y q)_p &= 0
\end{align}
where $J = \nabla_q f$, $K = \nabla_q g$ and $p$ is any point value at a node or edge midpoint, while $ (D_x q)_p $ is a finite-difference-type approximation to the derivative of $q$ at the location of the point value. If at $p$ a unique derivative $\del_x q_h$ or $\del_y q_h$ exists (tangential to the edge at the midpoint), then this is what is used; otherwise a(n upwind) regularization is used in complete analogy to the one-dimensional case \eqref{eq:pointvaluesAF}. One thus has:
\begin{align}
 \frac{\dd}{\dd t} q_{i+\frac12,j+\frac12} &+ J^+ (D_x^+ q)_{i+\frac12,j+\frac12} + J^- (D_x^- q)_{i+\frac12,j+\frac12} \\\nonumber&+ K^+ (D_y^+ q)_{i+\frac12,j+\frac12} + K^- (D_y^- q)_{i+\frac12,j+\frac12} = 0 \\
 \frac{\dd}{\dd t} q_{i+\frac12,j} &+ J^+ (D_x^+ q)_{i+\frac12,j} +J^- (D_x^- q)_{i+\frac12,j}  + K (D_y q)_{i+\frac12,j}  = 0 
\end{align}
and similarly for the perpendicular edge. It has been again proposed to obtain the finite differences by differentiating the reconstruction in the upwind cell:
\begin{align}
 (D_x^+ q)_{i+\frac12,j+\frac12} &= \frac{\del}{\del x} q_{\text{recon},ij}\left(\frac{\Delta x}{2}, \frac{\Delta y}{2}\right) = \frac{\del}{\del x} q_{\text{recon},i,j+1}\left(\frac{\Delta x}{2}, -\frac{\Delta y}{2}\right)\\
 (D_x^+ q)_{i+\frac12,j} &= \frac{\del}{\del x} q_{\text{recon},ij}\left(\frac{\Delta x}{2}, 0\right)\\
 (D_y q)_{i+\frac12,j} &= \frac{\del}{\del y} q_{\text{recon},ij}\left(\frac{\Delta x}{2}, 0\right) = \frac{\del}{\del y} q_{\text{recon},i+1,j}\left(-\frac{\Delta x}{2}, 0\right)
\end{align}
etc. The reader is referred to \cite{barsukow24afeuler,barsukow24affourier} for further details.

\subsection{Biorthogonal test functions and the Petrov-Galerkin method}

We seek biorthogonal test functions $\psi_{i+\frac12,j+\frac12}, \psi_{i+\frac12,j}  ,\psi_{i,j+\frac12}, \psi_{ij}$, such that
\begin{align}
(\psi_{r}, \phi_{s})_{L^2} &= \delta_{rs} \qquad \forall (r, s) \in \left\{ \left(i,j\right), \left(i+\frac12,0\right),\left (0,j+\frac12\right),\left (i+\frac12,j+\frac12\right)  \right\}
\end{align}

The space of test functions shall be denoted by
\begin{align}
 \mathscr W_\text{2d} := \mathrm{span}_{(i,j)\in\mathbb Z^2}\Big(\psi_{ij}, \psi_{i+\frac12,j}, \psi_{i,j+\frac12}, \psi_{i+\frac12, j+\frac12}\Big)
\end{align}

\begin{theorem}
 The test function $\psi_{ij}(x,y)$ associated to the average is $\frac{1}{\Delta x\Delta y} \id_{C_{ij}}$.
\end{theorem}
\begin{proof}
 Indeed,
 \begin{align}
  \int_{\mathbb R^2} \id_{C_{ij}} \phi_{i',j'} \dd \vec x &= \Delta x\Delta y \delta_{i,i'} \delta_{j,j'}\\
  \int_{\mathbb R^2} \id_{C_{ij}} \phi_{r} \dd \vec x &= 0  \qquad \forall r \in \left\{ \left(i'+\frac12,0\right),\left(0,j'+\frac12\right),\left(i'+\frac12,j'+\frac12\right) \right\}
 \end{align}
 by definition.
\end{proof}

Introduce again the local versions of the test functions, in complete analogy to the 1-d case, as 
\begin{align}
 A_0, A_{\pm\frac12,\pm\frac12}, A_{\pm\frac12,0}, A_{0,\pm\frac12} \colon \left[ -\frac{\Delta x}{2}, \frac{\Delta x}{2} \right ] \times \left[ -\frac{\Delta y}{2}, \frac{\Delta y}{2} \right ] \to \mathbb R
\end{align}

\begin{definition}
Given the numerical approximation \eqref{eq:numapprox2d}, the Active Flux Petrov-Galerkin method for the conservation law \eqref{eq:conslaw2d}, i.e. $\del_t q + \del_x f(q) + \del_y g(q) = 0$, is given by
\begin{align}
 (\psi, \del_t q_h)_{L^2} + (\psi, \del_x f(q_h))_{L^2}  +(\psi, \del_y g(q_h))_{L^2} = 0 \qquad \forall \psi \in \mathscr W_\text{2d}
\end{align}

\end{definition}

\subsubsection{Edge midpoints} \label{ssec:edgestests}

\begin{figure}
 \centering
 \includegraphics[width=0.7\textwidth]{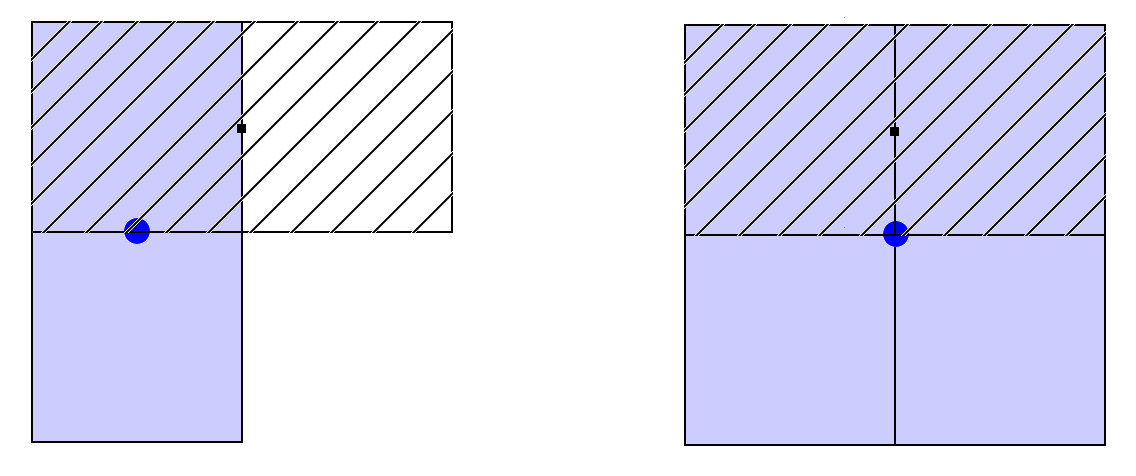}
 \caption{Illustration of intersections of supports between the test function $\psi_{i+\frac12,j}$ (diagonal lines) and different basis functions (shaded blue). \emph{Left}: One-cell intersection between supports of $\psi_{i+\frac12,j}$ and $\phi_{i,j-\frac12}$. The biorthogonality condition thus only involves the integral over $A_{+\frac12,0} B_{0,-\frac12}$. \emph{Right}: Two-cell intersection between supports of $\psi_{i+\frac12,j}$ and $\phi_{i+\frac12,j-\frac12}$. The biorthogonality condition involves the integrals over $A_{+\frac12,0} B_{+\frac12,-\frac12}$ and $A_{-\frac12,0} B_{-\frac12,-\frac12}$.}
 \label{fig:overlap-edge1}
\end{figure}

\begin{figure}
 \centering
 \includegraphics[width=0.35\textwidth]{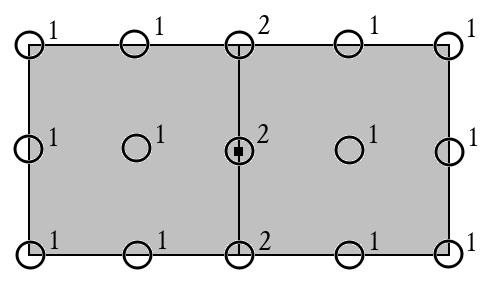}
 \caption{Numbers of cells in the intersection between supports of the test function $\psi_{i+\frac12,j}$ (center marked with a square) and various basis functions (centers marked with circles).}
 \label{fig:overlap-edge2}
\end{figure}

The test function $$\psi_{i+\frac12,j}(x,y) = A_{+\frac12,0}(x - x_i, y-y_j) \id_{C_{ij}} + A_{-\frac12,0}(x - x_{i+1},y - y_j) \id_{C_{i+1,j}}$$ associated to the point value at the edge midpoint $(i+\frac12,j)$ is assumed to be supported in $C_{ij} \cup C_{i+1,j}$. The biorthogonality conditions involve $L^2$ scalar products between the test function and various basis functions. These are all compactly supported, but the support varies between one cell (for the basis function associated to the cell average) and four cells (basis function associated to the node). It thus makes sense to organize the conditions according to the size of the intersection between the supports of the test and basis functions, starting with the basis functions whose supports intersect that of $\psi_{i+\frac12,j}$ only in one cell (see also Figures \ref{fig:overlap-edge1} and \ref{fig:overlap-edge2}):
\begin{align}
 \int_{C} A_{+\frac12,0} B_0 \,\dd \vec x &= 0 &
 \int_{C} A_{+\frac12,0} B_{-\frac12,0}  \,\dd \vec x &= 0&
 \int_{C} A_{+\frac12,0} B_{0,+\frac12}  \,\dd \vec x &= 0 \label{eq:2dAplus12first}\\
 \int_{C} A_{+\frac12,0} B_{0,-\frac12}  \,\dd \vec x &= 0&
 \int_{C} A_{+\frac12,0} B_{-\frac12,+\frac12}  \,\dd \vec x &= 0&
 \int_{C} A_{+\frac12,0} B_{-\frac12,-\frac12}  \,\dd \vec x &= 0 \label{eq:2dAplus12last}
\end{align}
where $C = [-\frac{\Delta x}{2}, \frac{\Delta x}{2}]\times [-\frac{\Delta y}{2}, \frac{\Delta y}{2}]$. Analogous formulae hold true for $A_{-\frac12,0}$. Additionally, there are three basis functions whose support includes all of the support of $\psi_{i+\frac12,j}$:
\begin{align}
 \int_{C} A_{+\frac12,0} B_{+\frac12,+\frac12}  \,\dd \vec x + \int_{C} A_{-\frac12,0} B_{-\frac12,+\frac12}  \,\dd \vec x &= 0 \label{eq:testfctedge1}\\
 \int_{C} A_{+\frac12,0} B_{+\frac12,-\frac12}  \,\dd \vec x + \int_{C} A_{-\frac12,0} B_{-\frac12,-\frac12}  \,\dd \vec x &= 0\\
 \int_{C} A_{+\frac12,0} B_{+\frac12,0}  \,\dd \vec x+\int_{C} A_{-\frac12,0} B_{-\frac12,0}  \,\dd \vec x &= 1 \label{eq:testfctedge3}
\end{align}

If $A_{\pm\frac12,0}$ are sought in $P^{2,2}$, 18 equations are required. These are only 15, there are thus 3 free parameters. Once again, one can therefore replace \eqref{eq:testfctedge1}--\eqref{eq:testfctedge3} by
\begin{align}
 \int_{C} A_{+\frac12,0} B_{+\frac12,+\frac12}  \,\dd \vec x &= \alpha_1 & \int_{C} A_{-\frac12,0} B_{-\frac12,+\frac12}  \,\dd \vec x &= -\alpha_1 \label{eq:2dAedge1} \\
 \int_{C} A_{+\frac12,0} B_{+\frac12,-\frac12}  \,\dd \vec x &= \alpha_2 &  \int_{C} A_{-\frac12,0} B_{-\frac12,-\frac12}  \,\dd \vec x &= -\alpha_2\\
 \int_{C} A_{+\frac12,0} B_{+\frac12,0}  \,\dd \vec x &= \frac12 + \frac{\alpha_3}{2} & \int_{C} A_{-\frac12,0} B_{-\frac12,0}  \,\dd \vec x &=  \frac12 - \frac{\alpha_3}{2} \label{eq:2dAedge1last}
\end{align}
with $\alpha_1,\alpha_2, \alpha_3 \in \mathbb R$. Quite generally, therefore, free parameters arise whenever the intersection between the support of the test function and that of the basis function is larger than one cell. Both $A_{+\frac12,0}$ and $A_{-\frac12,0}$ are now subject to nine equations each. The associated linear system has a unique solution. This is established as follows: The functions $\{ B_{\pm\frac12,\pm\frac12}, B_{\pm\frac12,0}, B_{0,\pm\frac12}, B_0\}$ form a basis of the space $P^{2,2}$, which can also be endowed with a tensor-product-basis $\{ L_s(x) L_r(y) \}_{s, r = 0,1,2}$ of Legendre polynomials. In particular, every element of the latter basis is a linear combination of elements of the former basis, such that the 9 equations in \eqref{eq:2dAplus12first}--\eqref{eq:2dAplus12last}, \eqref{eq:2dAedge1}--\eqref{eq:2dAedge1last} that involve $A_{+\frac12,0}$ can be written as
\begin{align}
 \int_{C} A_{+\frac12,0} L_s(x) L_r(y) \,\dd \vec x &= \mathrm{RHS} \qquad s,r = 0, 1, 2
\end{align}
But since $A_{+\frac12,0}$ is sought in the same space, it can also be expressed in the Legendre basis as
\begin{align}
 A_{+\frac12,0}(x, y) = \sum_{s',r' = 0}^2 a_{s'r'} L_{s'}(x) L_{r'}(y) 
\end{align}
Now, by orthogonality of the Legendre polynomials the linear system becomes diagonal with non-zero entries on the diagonal, since $C$ is a tensor product of two intervals:
\begin{align}
 \sum_{s',r' = 0}^2 a_{s'r'} \left(\int_{-\frac{\Delta x}{2}}^{\frac{\Delta x}{2}} L_{s'}(x)  L_s(x)\, \dd x \right ) \left( \int_{-\frac{\Delta y}{2}}^{\frac{\Delta y}{2}} L_{r'}(y)  L_r(y) \,\dd y \right) &= \mathrm{RHS} \\ \nonumber \qquad s,r = 0, 1, 2
\end{align}

\begin{figure}
 \centering
 \includegraphics[width=0.25\textwidth]{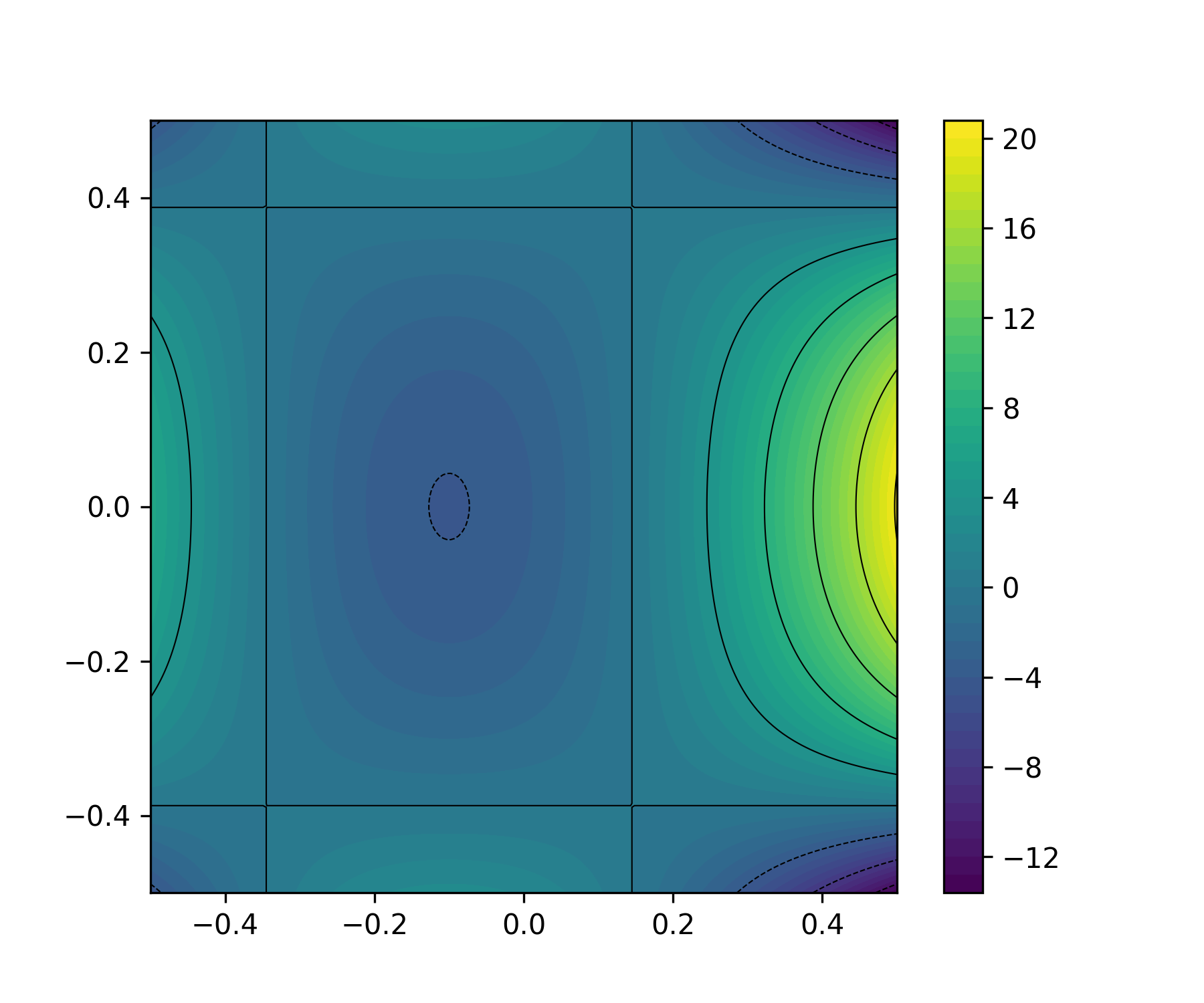}\hfill\hfill\includegraphics[width=0.25\textwidth]{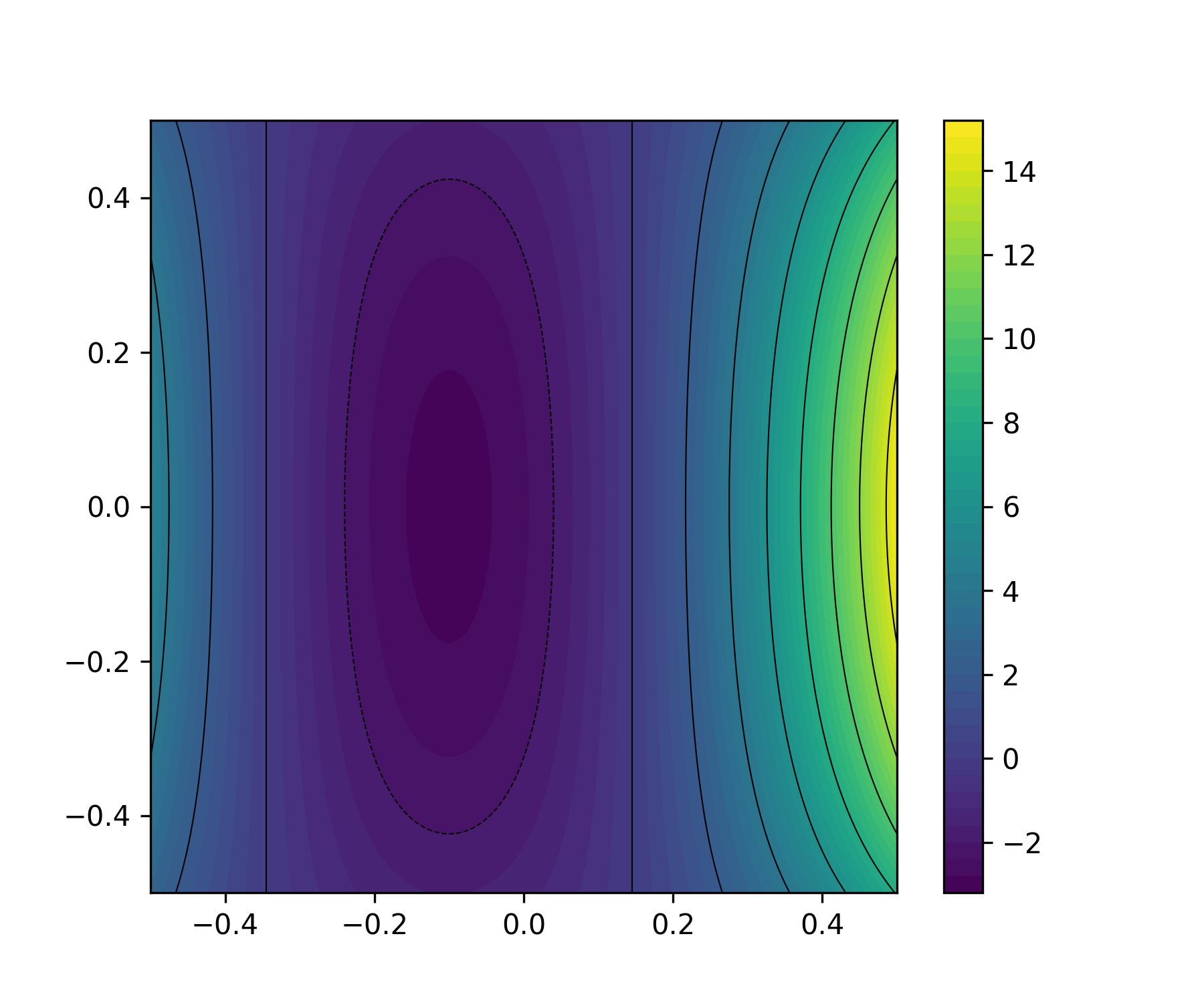}\hfill\includegraphics[width=0.25\textwidth]{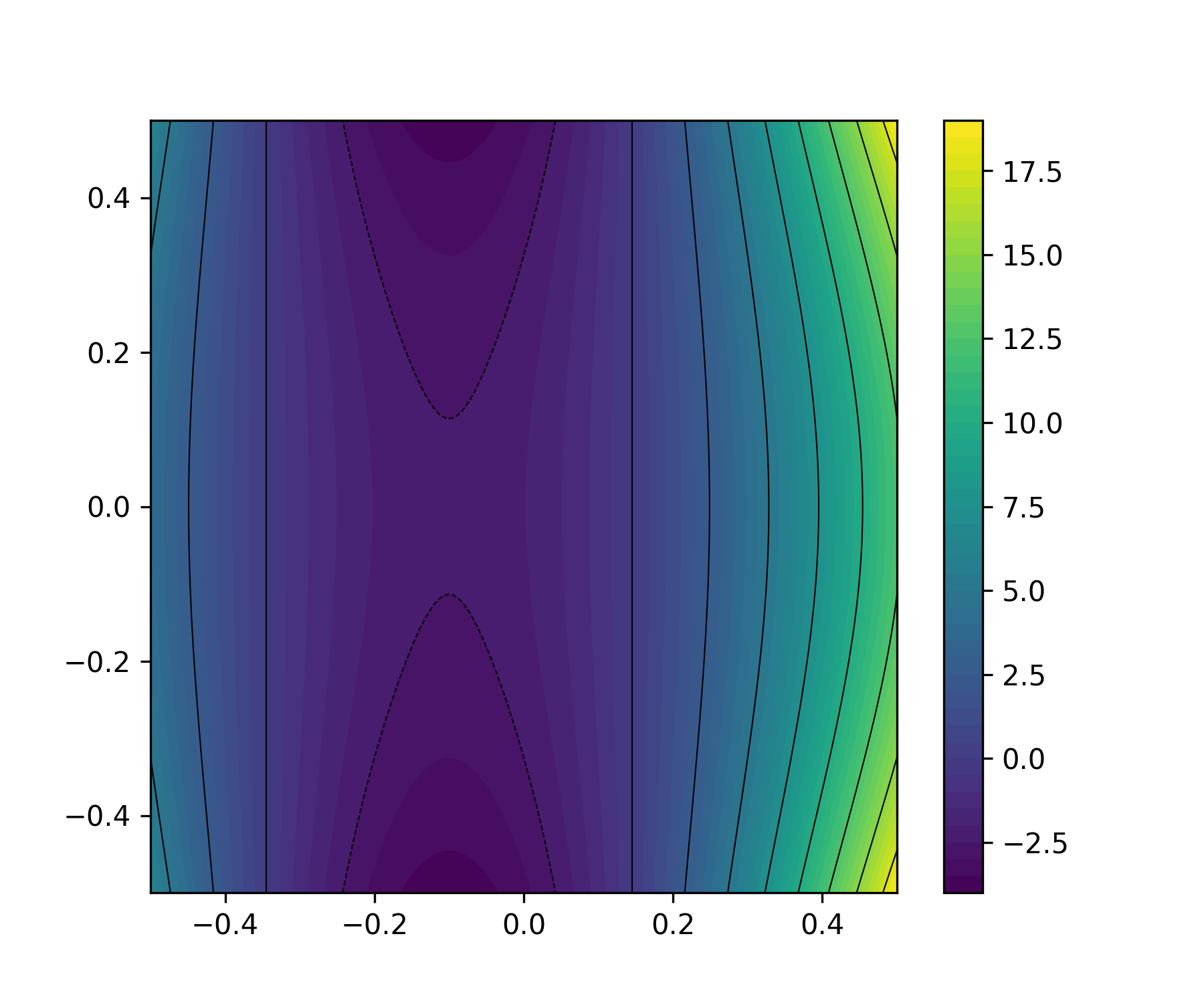}\hfill\includegraphics[width=0.25\textwidth]{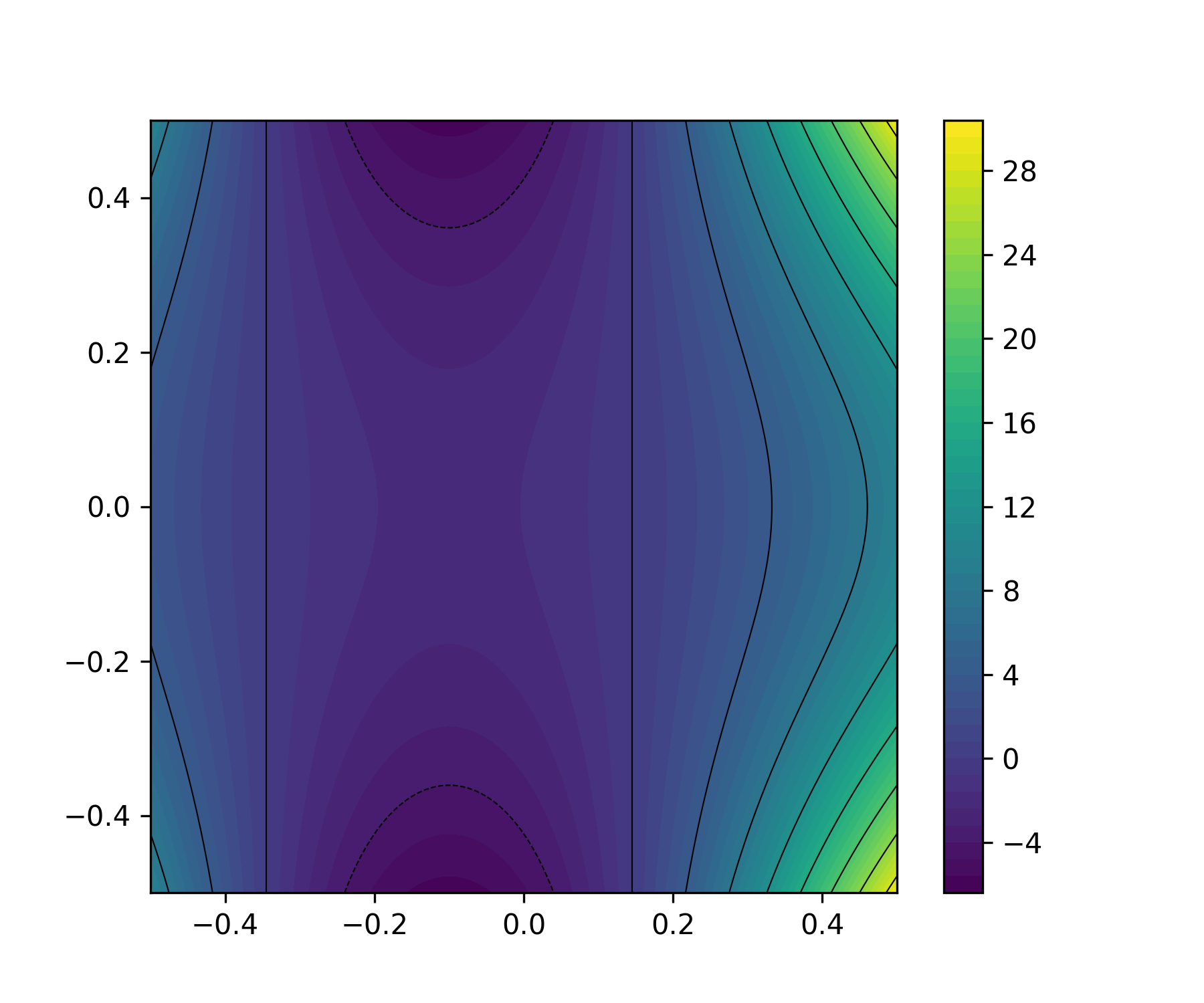} \\
 \includegraphics[width=0.25\textwidth]{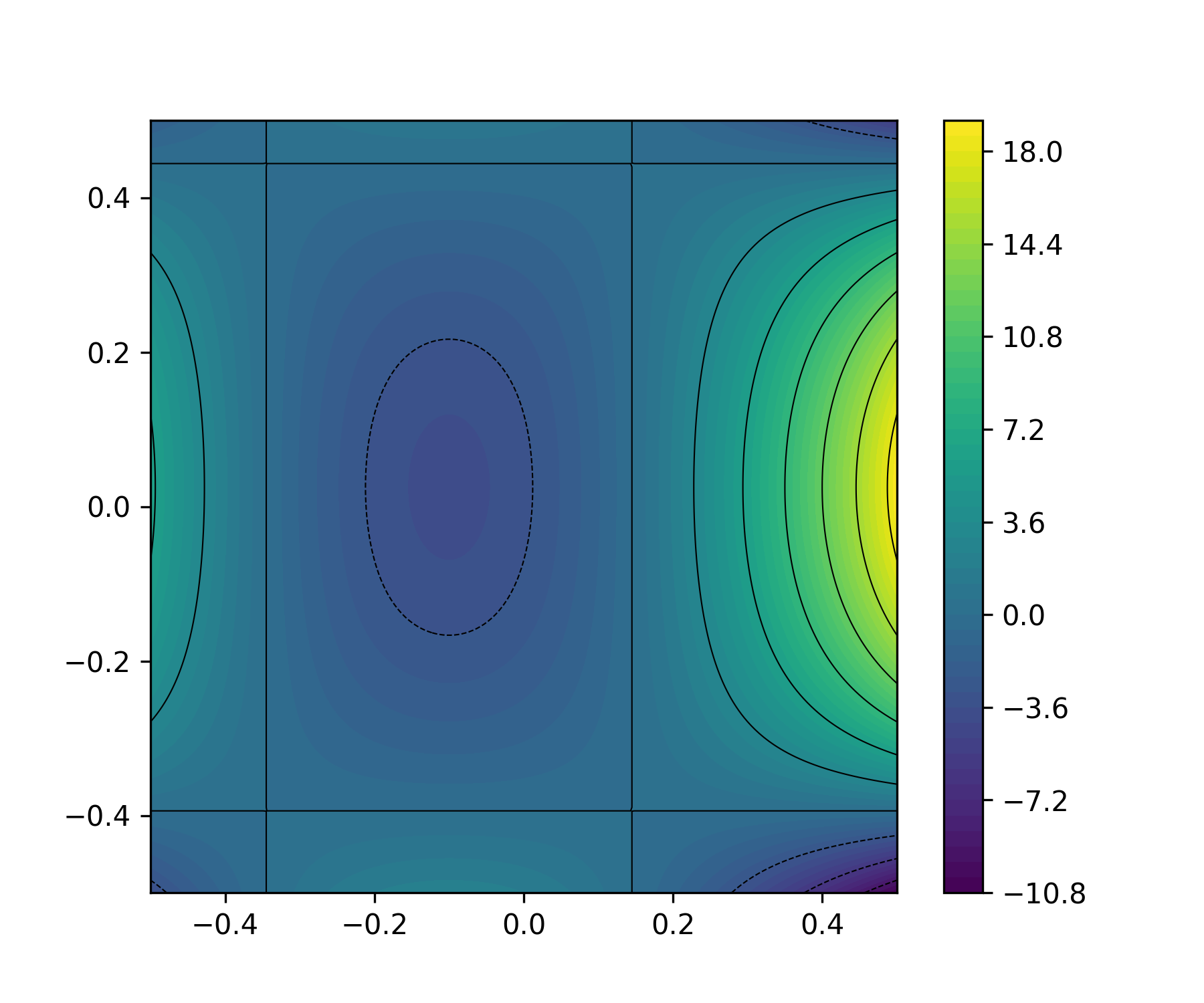}\hfill\hfill\includegraphics[width=0.25\textwidth]{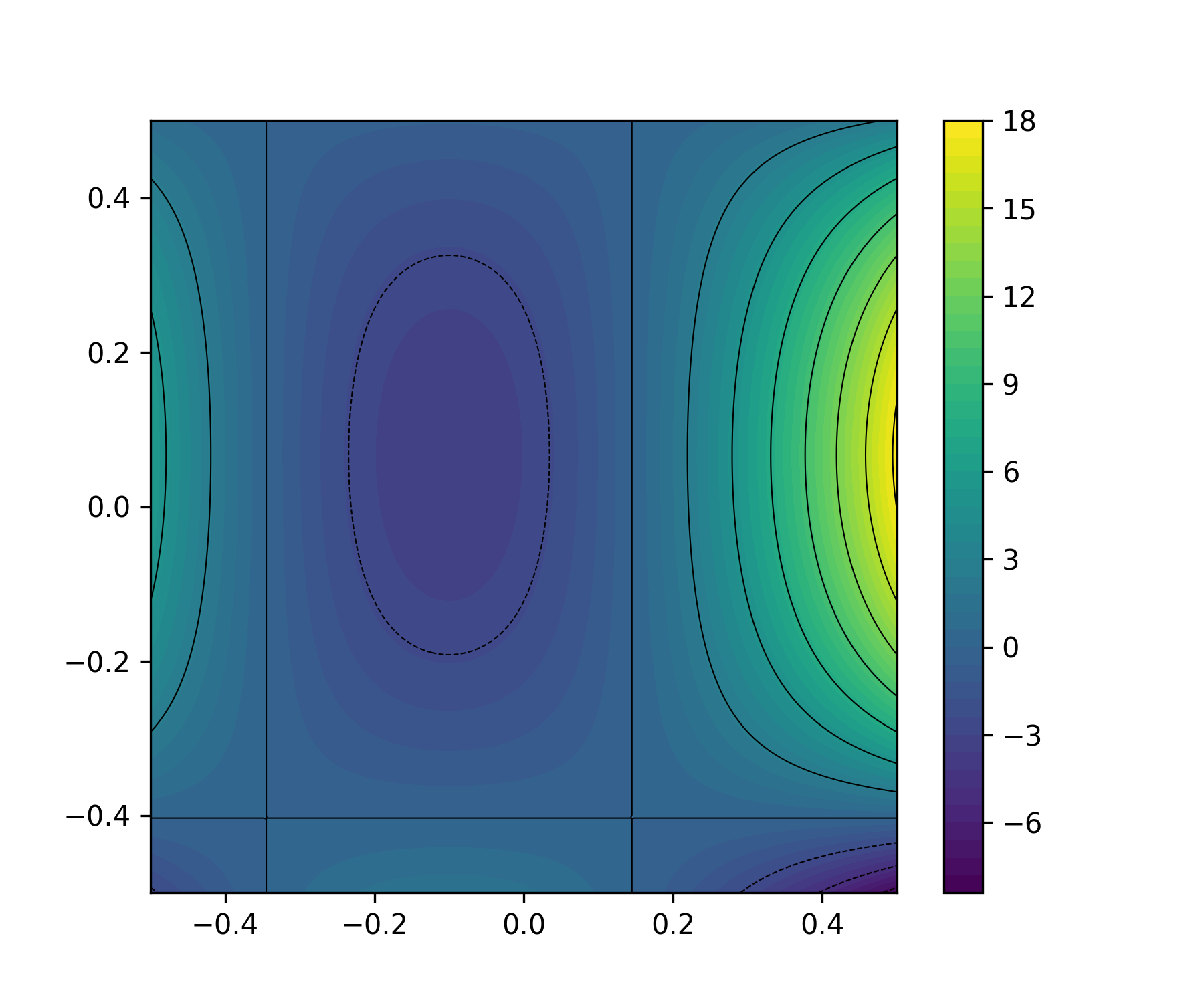}\hfill\includegraphics[width=0.25\textwidth]{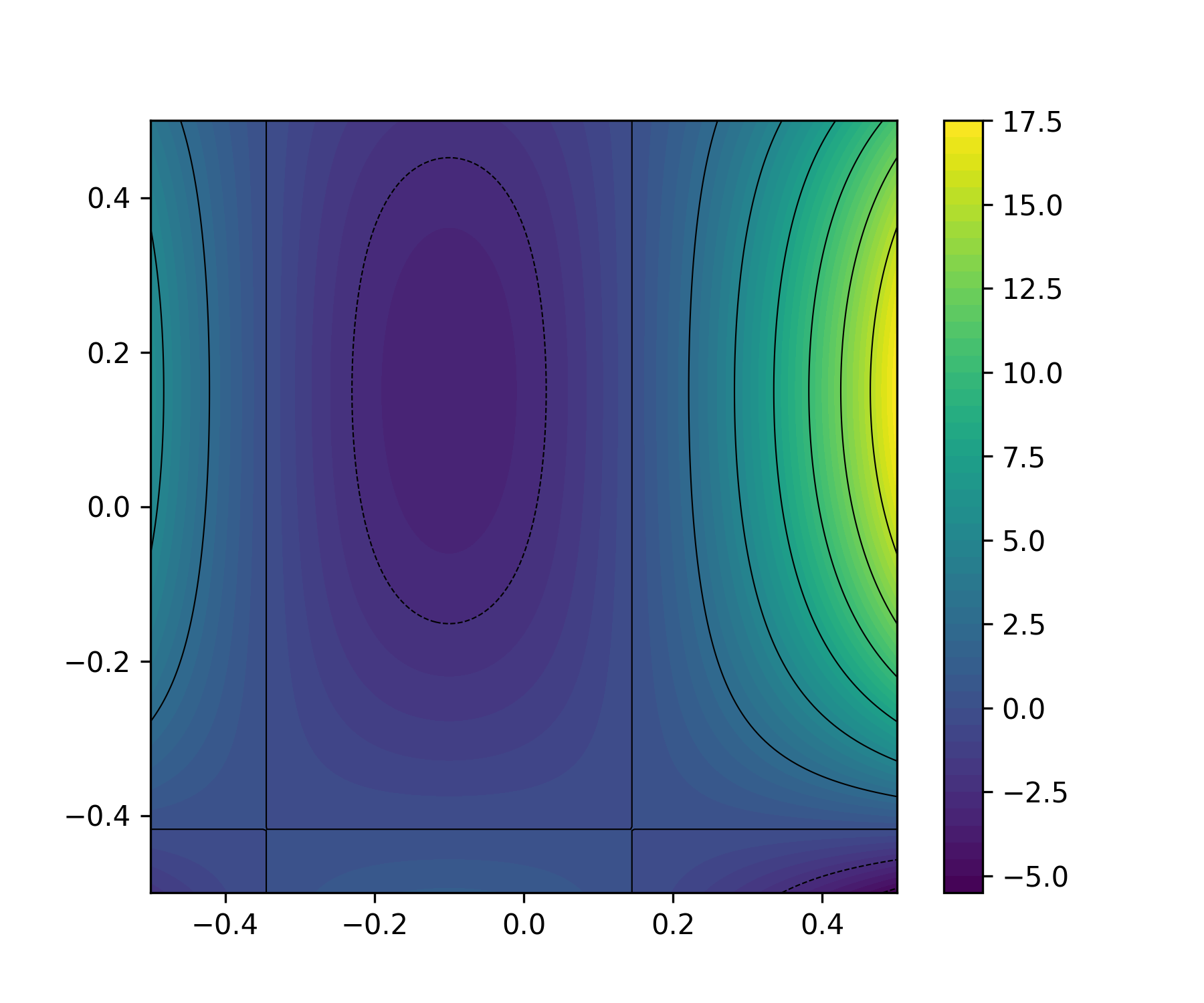}\hfill\includegraphics[width=0.25\textwidth]{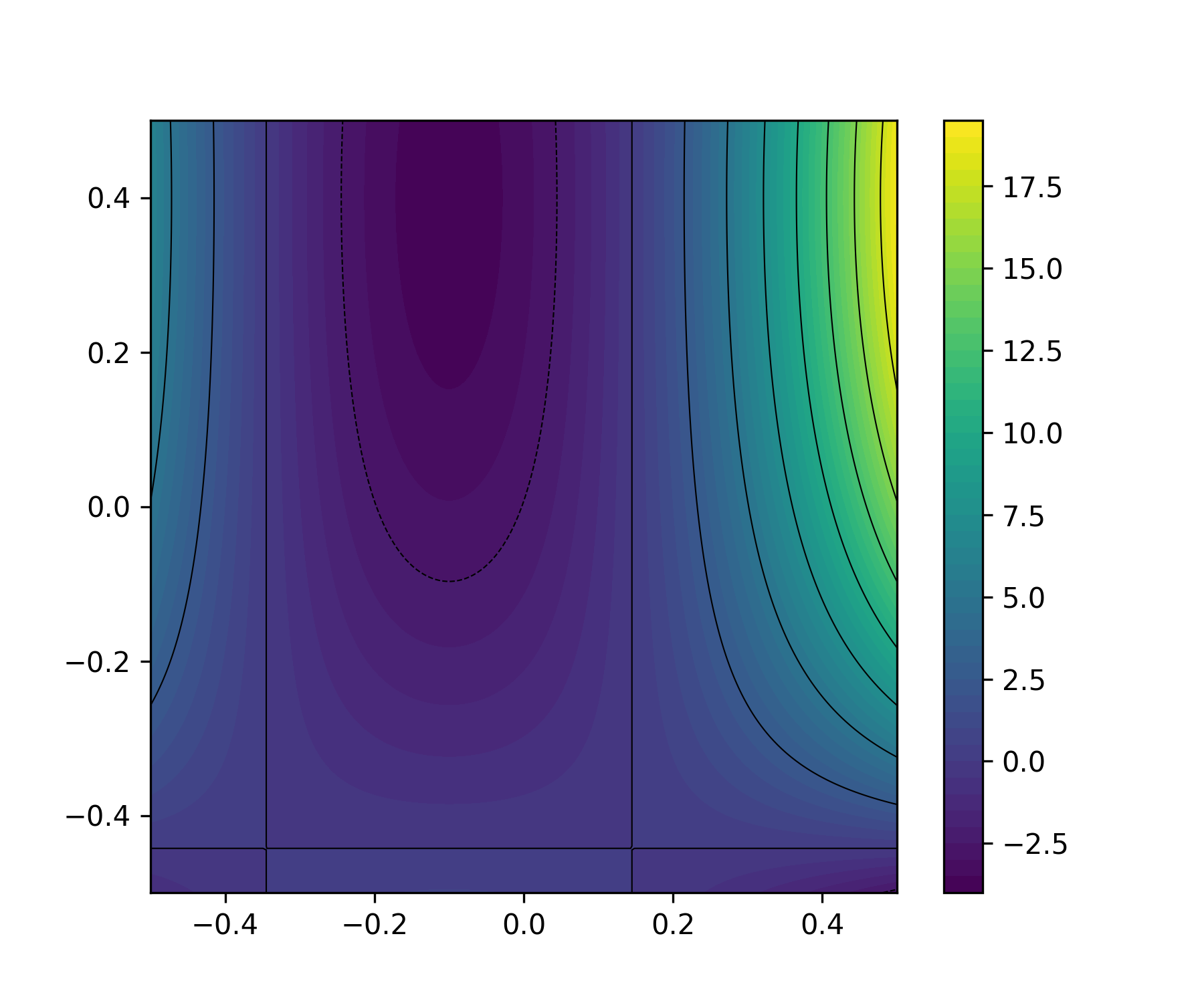}
 \caption{$A_{+\frac12,0}$ is shown for different parameter values, with $\alpha_3 = 1$ everywhere, i.e. fully upwinded. \emph{Top}: $\alpha_1 = \alpha_2 = 0, 0.2, 0.3, 0.4$ (\emph{from left to right}), $\alpha_3 = 1$. \emph{Bottom}: $\alpha_1 = 0.1, 0.2, 0.3, 0.4$ (\emph{from left to right}), $\alpha_2 = 0$.}
 \label{fig:2dtestfctedge}
\end{figure}

We refrain from solving the system now directly; instead, in Section \ref{ssec:tensorbasis} a more elegant way to explicitly obtain the test functions is shown. $A_{+\frac12,0}$ is shown in Figure \ref{fig:2dtestfctedge} for several values of $\alpha_1$, $\alpha_2$, $\alpha_3$. For reference, its explicit form is
\begin{align}
 A_{+\frac12,0} &= \left(-\frac34 + 3\xi + 15\xi^2 \right ) \cdot \left( \frac{3}{4} \big(- 4 (\alpha_{1} + \alpha_{2}) + 3 (\alpha_{3} + 1)\big) \right . \label{eq:Aplus120}  \\  \nonumber & \left.
 + 12 \eta  (\alpha_{1} - \alpha_{2})
 - 15 \eta^{2}  \big( -4 (\alpha_{1} + \alpha_{2} )+  \alpha_{3} + 1 \big ) \phantom{\frac{1}{1}} \!\!\!\!\! \right )
\end{align}
with $\xi = \frac{x}{\Delta x}$, $\eta = \frac{y}{\Delta y}$. The first bracket appears in the one-dimensional test function \eqref{eq:1dbiorthtestsymm}.

The situation is entirely analogous for the test functions $\psi_{i,j+\frac12}$ (i.e. $A_{0,+\frac12}$ and $A_{0,-\frac12}$) associated to the perpendicular edge midpoints. 
Before performing a similar analysis for the test functions associated to nodes, the evolution equation for the edge midpoints will be studied.

Similarly to the procedure in the Proof of Theorem \ref{thm:derivatives}, one can compute the normal derivative $(\psi_{i+\frac12,j}, \del_x q_h)_{L^2}$ by expanding $\del_x q_h$ in the interior of any cell in the Finite Element basis, e.g. as
\begin{align}
 \del_x q_h(x,y) \Big \vert _{C_{ij}} &= d^{(i,j)}_{i+\frac12,j+\frac12} B_{+\frac12,+\frac12}(x-x_i,y-y_j) + d^{(i,j)}_{i+\frac12,j} B_{+\frac12,0}(x-x_i,y-y_j) \\
 \nonumber & \!\!\!\!\!\!\!\!\!\!\!\!\!\!\!\!\!\!\!\!\!\!\!+ d^{(i,j)}_{i+\frac12,j-\frac12} B_{+\frac12,-\frac12}(x-x_i,y-y_j) + \text{further terms orthogonal to $A_{+\frac12,0}$}
\end{align}
Here,
\begin{align}
 d^{(i,j)}_{s,p} = \del_x q_h(x_s, y_p) \Big \vert _{C_{ij}}
\end{align}
is understood as the limit taken from inside cell $(i,j)$.
Using now the biorthogonality of the test functions one obtains
\begin{align}
 (&\psi_{i+\frac12,j}, \del_x q_h)_{L^2} = \int_C A_{+\frac12,0} \left( d^{(i,j)}_{i+\frac12,j+\frac12} B_{+\frac12,+\frac12} + d^{(i,j)}_{i+\frac12,j} B_{+\frac12,0} + d^{(i,j)}_{i+\frac12,j-\frac12} B_{+\frac12,-\frac12} \right ) \dd \vec x \\
 &\qquad\qquad\qquad + \nonumber \int_C A_{-\frac12,0} \left( d^{(i+1,j)}_{i+\frac12,j+\frac12} B_{-\frac12,+\frac12} + d^{(i+1,j)}_{i+\frac12,j} B_{-\frac12,0} + d^{(i+1,j)}_{i+\frac12,j-\frac12} B_{-\frac12,-\frac12} \right ) \dd \vec x \\
 &= 
  d^{(i,j)}_{i+\frac12,j} \left(\frac12 + \frac{\alpha_3}{2} \right)
   + d^{(i+1,j)}_{i+\frac12,j} \left(\frac12 - \frac{\alpha_3}{2} \right) \label{eq:xderivexpand}\\ &\nonumber \qquad \qquad 
   - \left(d^{(i+1,j)}_{i+\frac12,j+\frac12} - d^{(i,j)}_{i+\frac12,j+\frac12}\right) \alpha_1 
 - \left(d^{(i+1,j)}_{i+\frac12,j-\frac12}-d^{(i,j)}_{i+\frac12,j-\frac12}\right)\alpha_2 \\
 &\overset{\alpha_1 = \alpha_2 = 0}{=} \begin{cases}
               d^{(i,j)}_{i+\frac12,j} & \alpha_3 = 1 \\ d^{(i+1,j)}_{i+\frac12,j} & \alpha_3 = -1
              \end{cases}
\end{align}

Thus, if $\alpha_1 = \alpha_2 = 0$, the discrete derivative arising in the Petrov-Galerkin framework is a combination of the derivatives $d^{(i,j)}_{i+\frac12,j}$ of the approximation from the left and $d^{(i+1,j)}_{i+\frac12,j}$ (from the right), according to the value of $\alpha_3$. The derivative is taken purely from the left (right) if $\alpha_3 = 1$ ($=-1$). The upwinding considered so far for Active Flux (e.g. in \cite{barsukow24afeuler}) for the point value at an edge was a choice between the derivative from one or the other side of the edge at the edge midpoint. As is customary, this upwind derivative can be expressed as the average of left and right derivatives minus a term proportional to the jump. $\alpha_3$ is just the usual parameter governing the amount of upwinding/numerical diffusion.

One can give a similar meaning to the free parameters $\alpha_1, \alpha_2$, even though choices $\alpha_1 \neq 0$, or $\alpha_2 \neq 0$ have so far not been considered in the literature. These parameters premultiply jumps in the same direction (normal to the edge), but evaluated at locations along the edge other than the midpoint (at the nodes, in fact). They are structurally similar to jump terms that contribute to stability.     
In practice, so far these terms do not seem to make a significant impact; a more detailed analysis of their role is subject of future work, but an example is provided in Section \ref{sec:numerical}. Analogous results follow for the other normal derivative $(\psi_{i,j+\frac12}, \del_y q_h)_{L^2}$.

The jumps of the derivatives of the approximation vanish if the solution is globally polynomial of degree at most $K$, i.e. the parameters $\alpha_1$, $\alpha_2$, $\alpha_3$ do not change the order of accuracy of the method.

For the tangential derivative $(\psi_{i+\frac12,j}, \del_y q_h)_{L^2}$, it has already been observed in the Active Flux literature that uniqueness of the trace of the approximations in the two adjacent cells means that there is no possibility to introduce upwinding. In other words, at the edge midpoint, the derivatives tangential to the edge are central ones. Define first the derivatives as limits from inside the cell $(i,j)$ as
\begin{align}
 t^{(i,j)}_{s,p} := \del_y q_h(x_s, y_p) \Big \vert _{C_{ij}} \label{eq:defyderiv}
\end{align}

Indeed, then first in complete analogy to \eqref{eq:xderivexpand} one finds again
\begin{align}
(\psi_{i+\frac12,j}, \del_y q_h)_{L^2} &= t^{(i,j)}_{i+\frac12,j} \left(\frac12 + \frac{\alpha_3}{2} \right)
   + t^{(i+1,j)}_{i+\frac12,j} \left(\frac12 - \frac{\alpha_3}{2} \right) \\ &\nonumber \qquad 
   - \left(t^{(i+1,j)}_{i+\frac12,j+\frac12} - t^{(i,j)}_{i+\frac12,j+\frac12}\right) \alpha_1 
 - \left(t^{(i+1,j)}_{i+\frac12,j-\frac12}-t^{(i,j)}_{i+\frac12,j-\frac12}\right)\alpha_2
 \end{align}
but now continuity of derivatives across the edge implies simply
\begin{align}
(\psi_{i+\frac12,j}, \del_y q_h)_{L^2} &= t^{(i,j)}_{i+\frac12,j}  = t^{(i+1,j)}_{i+\frac12,j}
\end{align}
This is the prescription used in the Active Flux literature (see e.g. \cite{barsukow24afeuler}).

\subsubsection{Nodes}

Similarly, we assume $\psi_{i+\frac12,j+\frac12}$ to have support in $C_{ij} \cup C_{i+1,j}  \cup C_{i,j+1} \cup C_{i+1,j+1}$:
\begin{align}
 \psi_{i+\frac12,j+\frac12}(x,y) &= A_{+\frac12,+\frac12}(x - x_i, y-y_j) \id_{C_{ij}} + A_{-\frac12,+\frac12}(x - x_{i+1},y - y_j) \id_{C_{i+1,j}} \\
 & \nonumber \!\!\!\!\!\!\!\!\!\!\!\!\!\!\!\!\!\!\!\!\!\!\!\!\!\! + A_{+\frac12,-\frac12}(x - x_i, y-y_{j+1}) \id_{C_{i,j+1}} + A_{-\frac12,-\frac12}(x - x_{i+1},y - y_{j+1}) \id_{C_{i+1,j+1}}
\end{align}

\begin{figure}
 \centering
 \includegraphics[width=0.7\textwidth]{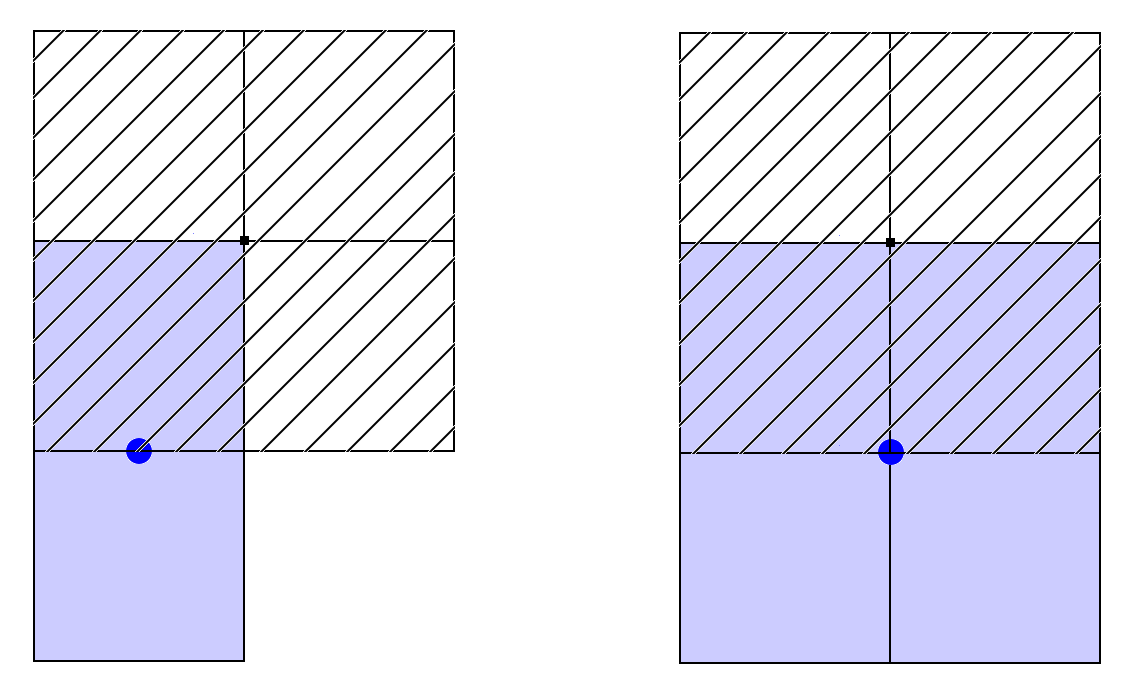}
 \caption{Illustration of intersections of supports between the test function $\psi_{i+\frac12,j+\frac12}$ (diagonal lines) and different basis functions (shaded blue). \emph{Left}: One-cell intersection between supports of $\psi_{i+\frac12,j+\frac12}$ and $\phi_{i,j-\frac12}$. The biorthogonality condition thus only involves the integral over $A_{+\frac12,+\frac12} B_{0,-\frac12}$. \emph{Right}: Two-cell intersection between supports of $\psi_{i+\frac12,j+\frac12}$ and $\phi_{i+\frac12,j-\frac12}$. The biorthogonality condition involves the integrals over $A_{+\frac12,+\frac12} B_{+\frac12,-\frac12}$ and $A_{-\frac12,+\frac12} B_{-\frac12,-\frac12}$. The four-cell intersection (omitted) is between the supports of $\psi_{i+\frac12,j+\frac12}$ and $\phi_{i+\frac12,j+\frac12}$ only.}
 \label{fig:overlap-node1}
\end{figure}

\begin{figure}
 \centering
 \includegraphics[width=0.35\textwidth]{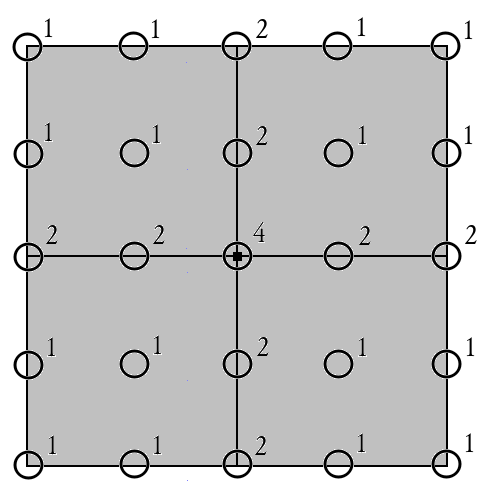}
 \caption{Numbers of cells in the intersection between supports of the test function $\psi_{i+\frac12,j+\frac12}$ (center marked with a square) and various basis functions (centers marked with circles).}
 \label{fig:overlap-node2}
\end{figure}

Again, the biorthogonality conditions are organized according to the size of the intersection between the supports of the test function and the basis function. Starting with basis functions which lead to an intersection of supports in only one cell, one has the following defining relations for the test functions $A_{\pm\frac12,\pm\frac12}$ associated to a node (see also Figures \ref{fig:overlap-node1} and \ref{fig:overlap-node2}):
\begin{align}
 \int_{C} A_{+\frac12,+\frac12} B_0 \,\dd \vec x &= 0 \label{eq:nodal1support1}&
 \int_{C} A_{+\frac12,+\frac12} B_{-\frac12,0}  \,\dd \vec x &= 0\\
 \int_{C} A_{+\frac12,+\frac12} B_{0,-\frac12}  \,\dd \vec x &= 0&
 \int_{C} A_{+\frac12,+\frac12} B_{-\frac12,-\frac12}  \,\dd \vec x &= 0 \label{eq:nodal1support4}
\end{align}
and similarly for the other three test functions, i.e. a total of 16 equations.
Two cells are involved in the following:
\begin{align}
 \int_{C} A_{+\frac12,+\frac12} B_{+\frac12,0}  \,\dd \vec x+ \int_{C} A_{-\frac12,+\frac12} B_{-\frac12,0}  \,\dd \vec x &= 0\\
 \int_{C} A_{+\frac12,+\frac12} B_{0,\frac12}  \,\dd \vec x + \int_{C} A_{+\frac12,-\frac12} B_{0,-\frac12}  \,\dd \vec x &= 0\\
 \int_{C} A_{+\frac12,+\frac12} B_{+\frac12,-\frac12}  \,\dd \vec x + \int_{C} A_{-\frac12,+\frac12} B_{-\frac12,-\frac12}  \,\dd \vec x &= 0\\
 \int_{C} A_{+\frac12,+\frac12} B_{-\frac12,+\frac12}  \,\dd \vec x + \int_{C} A_{+\frac12,-\frac12} B_{-\frac12,-\frac12}  \,\dd \vec x &= 0
\end{align}
Each of the relations exists analogously for another pair of test functions, i.e. a total of 8 equations. The remaining equation is
\begin{align}
 \int_{C} A_{+\frac12,+\frac12} B_{+\frac12,+\frac12}  \,\dd \vec x  &+ \int_{C} A_{-\frac12,+\frac12} B_{-\frac12,+\frac12}  \,\dd \vec x \label{eq:nodalself}\\
 \nonumber &+ \int_{C} A_{+\frac12,-\frac12} B_{+\frac12,-\frac12}  \,\dd \vec x + \int_{C} A_{-\frac12,-\frac12} B_{-\frac12,-\frac12}  \,\dd \vec x = 1 
\end{align}

These are 25 equations for the $4 \times 9$ parameters of the four test functions, leaving 11 of them free. For $A_{+\frac12,+\frac12} $ one thus again introduces the following equations:
\begin{align}
 \int_{C} A_{+\frac12,+\frac12} B_{+\frac12,0}  \,\dd \vec x &= \alpha_1\\
 \int_{C} A_{+\frac12,+\frac12} B_{0,\frac12}  \,\dd \vec x &= \alpha_2\\
 \int_{C} A_{+\frac12,+\frac12} B_{+\frac12,-\frac12}  \,\dd \vec x &= \alpha_3\\
 \int_{C} A_{+\frac12,+\frac12} B_{-\frac12,+\frac12}  \,\dd \vec x &= \alpha_4
\end{align}
Each of the relations exists analogously for another pair of test functions:
\begin{align}
 \int_{C} A_{-\frac12,+\frac12} B_{0,+\frac12} \,\dd \vec x = \alpha_5 = - \int_{C} A_{-\frac12,-\frac12} B_{0,-\frac12} \,\dd \vec x \\
 \int_{C} A_{+\frac12,-\frac12} B_{+\frac12,0} \,\dd \vec x = \alpha_6 = - \int_{C} A_{-\frac12,-\frac12} B_{-\frac12,0} \,\dd \vec x \\
 \int_{C} A_{-\frac12,+\frac12} B_{+\frac12,+\frac12} \,\dd \vec x = \alpha_7 = - \int_{C} A_{-\frac12,-\frac12} B_{+\frac12,-\frac12} \,\dd \vec x \\
 \int_{C} A_{+\frac12,-\frac12} B_{+\frac12,+\frac12} \,\dd \vec x = \alpha_8 = - \int_{C} A_{-\frac12,-\frac12} B_{-\frac12,+\frac12} \,\dd \vec x 
\end{align}

These are a total of 8 equations. The remaining equation is
\begin{align}
 \int_{C} A_{+\frac12,+\frac12} B_{+\frac12,+\frac12}  \,\dd \vec x &= \frac14 + \frac{\alpha_9}{4} + \alpha_{10}
\end{align}
having split \eqref{eq:nodalself} as
\begin{align}
 \int_{C} A_{+\frac12,+\frac12} B_{+\frac12,+\frac12}  \,\dd \vec x  + \int_{C} A_{-\frac12,+\frac12} B_{-\frac12,+\frac12} \,\dd \vec x  &= \frac12 + \frac{\alpha_9}{2} \\
 \int_{C} A_{+\frac12,-\frac12} B_{+\frac12,-\frac12}  \,\dd \vec x+ \int_{C} A_{-\frac12,-\frac12} B_{-\frac12,-\frac12}  \,\dd \vec x &= \frac12 -\frac{\alpha_9}{2} 
\end{align}

Observe again that free parameters are linked to the support of the test function intersecting that of a basis function in more than one cell.
Together with \eqref{eq:nodal1support1}--\eqref{eq:nodal1support4} they give a unique solution, since every one of the four constituent polynomials $A_{\pm\frac12,\pm\frac12} \in P^{2,2}$ is subject to nine equations of the form
\begin{align}
 \int_C A_{\pm\frac12,\pm\frac12} B_{k,\ell} \dd \vec x &= \mathrm{RHS}.
\end{align}
Uniqueness then follows by the same argument as the one given in Section \ref{ssec:edgestests} for the test functions associated to the edges.

$(\psi_{i+\frac12,j+\frac12}, \del_x q_h)_{L^2}$ can now also be interpreted as a weighted average of the derivatives of $q_h$ in various points:

\begin{align}
 &(\psi_{i+\frac12,j+\frac12}, \del_x q_h)_{L^2}  \\ \nonumber &= \int_C A_{+\frac12,+\frac12} \Big( 
 d^{(i,j)}_{i-\frac12,j+\frac12} B_{-\frac12,+\frac12}
 +d^{(i,j)}_{i,j+\frac12} B_{0,+\frac12}
 +d^{(i,j)}_{i+\frac12,j+\frac12} B_{+\frac12,+\frac12}
 \\& \phantom{mmmmmmmmmm} \nonumber 
 +d^{(i,j)}_{i+\frac12,j} B_{+\frac12,0}
 +d^{(i,j)}_{i+\frac12,j-\frac12} B_{+\frac12,-\frac12}
 \Big) \dd \vec x\\ \nonumber
 &+\int_C A_{-\frac12,+\frac12} \Big( 
 d^{(i+1,j)}_{i+\frac32,j+\frac12} B_{+\frac12,+\frac12}
 +d^{(i+1,j)}_{i+1,j+\frac12} B_{0,+\frac12}
 +d^{(i+1,j)}_{i+\frac12,j+\frac12} B_{-\frac12,+\frac12}
 \\& \phantom{mmmmmmmmmm} \nonumber 
 +d^{(i+1,j)}_{i+\frac12,j} B_{-\frac12,0}
 +d^{(i+1,j)}_{i+\frac12,j-\frac12} B_{-\frac12,-\frac12}
 \Big) \dd \vec x\\ \nonumber
 &+ \int_C A_{+\frac12,-\frac12} \Big( 
 d^{(i,j+1)}_{i-\frac12,j+\frac12} B_{-\frac12,-\frac12}
 +d^{(i,j+1)}_{i,j+\frac12} B_{0,-\frac12}
 +d^{(i,j+1)}_{i+\frac12,j+\frac12} B_{+\frac12,-\frac12}
 \\& \phantom{mmmmmmmmmm} \nonumber 
 +d^{(i,j+1)}_{i+\frac12,j+1} B_{+\frac12,0}
 +d^{(i,j+1)}_{i+\frac12,j+\frac32} B_{+\frac12,+\frac12}
 \Big) \dd \vec x \\ \nonumber
 &+ \int_C A_{-\frac12,-\frac12} \Big( 
 d^{(i+1,j+1)}_{i+\frac32,j+\frac12} B_{+\frac12,-\frac12}
 +d^{(i+1,j+1)}_{i+1,j+\frac12} B_{0,-\frac12}
 +d^{(i+1,j+1)}_{i+\frac12,j+\frac12} B_{-\frac12,-\frac12}
 \\& \phantom{mmmmmmmmmm} \nonumber 
 +d^{(i+1,j+1)}_{i+\frac12,j+1} B_{-\frac12,0}
 +d^{(i+1,j+1)}_{i+\frac12,j+\frac32} B_{-\frac12,+\frac12}
 \Big) \dd \vec x \\ 
 &
 =  
 d^{(i,j)}_{i+\frac12,j+\frac12} \left(\frac14 + \frac{\alpha_9}{4} + \alpha_{10}\right)
 +d^{(i+1,j)}_{i+\frac12,j+\frac12} \left( \frac14 + \frac{\alpha_9}{4} - \alpha_{10} \right) 
 \\&\!\!\!\!\!\!\!\!\!\!\!\!\!\!\!\!\!\!\!\! \nonumber\qquad  +d^{(i,j+1)}_{i+\frac12,j+\frac12} \left(\frac14 - \frac{\alpha_9}{4} + \alpha_{11}\right)
 +d^{(i+1,j+1)}_{i+\frac12,j+\frac12} \left( \frac14 - \frac{\alpha_9}{4} - \alpha_{11}   \right)
 \\&\!\!\!\!\!\!\!\!\!\!\!\!\!\!\!\!\!\!\!\! \nonumber\qquad +\alpha_1 \left(d^{(i,j)}_{i+\frac12,j}  -d^{(i+1,j)}_{i+\frac12,j} \right)
 +\alpha_2 \left( d^{(i,j)}_{i,j+\frac12}  -d^{(i,j+1)}_{i,j+\frac12}\right)
 +\alpha_3 \left( d^{(i,j)}_{i+\frac12,j-\frac12}  -d^{(i+1,j)}_{i+\frac12,j-\frac12} \right)
 \\&\!\!\!\!\!\!\!\!\!\!\!\!\!\!\!\!\!\!\!\! \nonumber\qquad + \alpha_4 \left( d^{(i,j)}_{i-\frac12,j+\frac12} -d^{(i,j+1)}_{i-\frac12,j+\frac12} \right)
 +\alpha_5\left(d^{(i+1,j)}_{i+1,j+\frac12} -d^{(i+1,j+1)}_{i+1,j+\frac12} \right)
 \\&\!\!\!\!\!\!\!\!\!\!\!\!\!\!\!\!\!\!\!\! \nonumber\qquad +\alpha_6 \left( d^{(i,j+1)}_{i+\frac12,j+1} -d^{(i+1,j+1)}_{i+\frac12,j+1} \right)
 +\alpha_7\left( d^{(i+1,j)}_{i+\frac32,j+\frac12}  -d^{(i+1,j+1)}_{i+\frac32,j+\frac12} \right)
 \\&\!\!\!\!\!\!\!\!\!\!\!\!\!\!\!\!\!\!\!\! \nonumber\qquad +\alpha_8 \left(d^{(i,j+1)}_{i+\frac12,j+\frac32} -d^{(i+1,j+1)}_{i+\frac12,j+\frac32} \right) 
 \end{align}

 Observe that the expression $(\psi_{i+\frac12,j+\frac12}, \del_y q_h)_{L^2}$ has exactly the same form, with $d$ to be replaced by $t$ from \eqref{eq:defyderiv}. But now, due to global continuity of $q_h$ certain derivatives have no jump (e.g. $d^{(i,j)}_{i,j+\frac12}  -d^{(i,j+1)}_{i,j+\frac12}=0$), and this is where the two expressions become different:
 \begin{align}
&(\psi_{i+\frac12,j+\frac12}, \del_x q_h)_{L^2}
 = 
 d^{(i,j)}_{i+\frac12,j+\frac12} \left(\frac12  + \alpha_{10} + \alpha_{11}\right)
 +d^{(i+1,j)}_{i+\frac12,j+\frac12} \left( \frac12  - (\alpha_{10} + \alpha_{11})   \right) \label{eq:nodederivx}
 \\& \nonumber\qquad 
 +\alpha_1 \left(d^{(i,j)}_{i+\frac12,j}  -d^{(i+1,j)}_{i+\frac12,j} \right)
 +\alpha_3 \left( d^{(i,j)}_{i+\frac12,j-\frac12}  -d^{(i+1,j)}_{i+\frac12,j-\frac12} \right)
 \\& \nonumber\qquad +\alpha_6 \left( d^{(i,j+1)}_{i+\frac12,j+1} -d^{(i+1,j+1)}_{i+\frac12,j+1} \right)
 +\alpha_8 \left(d^{(i,j+1)}_{i+\frac12,j+\frac32} -d^{(i+1,j+1)}_{i+\frac12,j+\frac32} \right)
\end{align}
The parameters that have disappeared here are present for the derivative in $y$-direction:
\begin{align}
&(\psi_{i+\frac12,j+\frac12}, \del_y q_h)_{L^2}
=  
 t^{(i,j)}_{i+\frac12,j+\frac12} \left(\frac12 + \frac{\alpha_9}{2}  \right)  +t^{(i,j+1)}_{i+\frac12,j+\frac12} \left(\frac12 - \frac{\alpha_9}{2}    \right) \label{eq:nodederivy}
 \\& \nonumber\qquad 
 +\alpha_2 \left( t^{(i,j)}_{i,j+\frac12}  -t^{(i,j+1)}_{i,j+\frac12}\right)
 + \alpha_4 \left( t^{(i,j)}_{i-\frac12,j+\frac12} -t^{(i,j+1)}_{i-\frac12,j+\frac12} \right)
 \\& \nonumber\qquad 
 +\alpha_5\left(t^{(i+1,j)}_{i+1,j+\frac12} -t^{(i+1,j+1)}_{i+1,j+\frac12} \right)
 +\alpha_7\left( t^{(i+1,j)}_{i+\frac32,j+\frac12}  -t^{(i+1,j+1)}_{i+\frac32,j+\frac12} \right)
 \end{align}

Each discrete derivative contains 5 independent free parameters that control upwinding at the node $(i+\frac12,j+\frac12)$ and at the midpoints and nodes along $x_{i+\frac12}$/$y_{j+\frac12}$. So far, Active Flux methods in the literature only considered the upwinding at the node itself (i.e. the choice $\alpha_1 = \alpha_2 = \dots = \alpha_7 = \alpha_8 = 0$). A detailed study of the consequences of including the additional jumps is subject of future work.

\subsection{Relation to the tensor basis} \label{ssec:tensorbasis}

This section shows a simple way to derive explicit expressions for the test functions in 2-d.
Observe that the basis functions \eqref{eq:basis2d1}--\eqref{eq:basis2dlast} do not form a tensor basis. For example, (with $\xi := x/\Delta x$, $\eta := y/\Delta y$),
\begin{align}
 B_{+\frac12,0}(x,y) &= \frac14(2\xi+1)(6\xi-1)(1-4\eta^2) \\ &\neq \frac14 (2 \xi+1)(6\xi-1) \cdot \frac32 (1-4\eta^2)  = B_{+\frac12}(x) B_0(y) \\
 B_{+\frac12,+\frac12}(x,y) &= \frac1{16} (2\xi+1)(2\eta+1)(12 \xi \eta + 2 \xi + 2 \eta -1) \\ &\neq \frac1{16} (2\xi+1)(2\eta+1) (6\xi-1)(6\eta-1) = B_{+\frac12}(x)B_{+\frac12}(y)
\end{align}
Also the test functions in 2-d are not tensor-product expressions of 1-d test functions.
However, it is easy to find the expressions of the test functions in the tensor basis. 

Consider the tensor basis $\tilde B_{r,s}(x,y) := B_r(x) B_s(y)$ of $Q^2$. In fact, $ B_{+\frac12,0}(x,y)$ is not $B_{+\frac12}(x) B_0(y)$ because $B_0(0) \neq 1$. Instead, we have

\begin{theorem}
 The 2-d basis functions $B$ of Active Flux can be expressed in the tensor basis $\tilde B$ as follows:
 \begin{align}
 B_{r,s} &= \frac23 \tilde B_{r,s} && \forall (r,s) \in \left\{ \left(\pm\frac12,0\right), \left(0, \pm\frac12\right) \right \} && \text{(edges)} \label{eq:tensorbasisedge} \\
 B_{r,s} &= \tilde B_{r,s} + \frac16 \tilde B_{r,0} + \frac16 \tilde B_{0,s} && \forall (r,s) \in \left\{ \left(\pm\frac12,\pm\frac12\right) \right \} && \text{(nodes)} \label{eq:tensorbasisnode} \\
 B_{0} &= \tilde B_{0,0} \label{eq:tensorbasisavg}
\end{align}
\end{theorem}
\begin{proof}
Write $x_{\pm\frac12} = \pm\frac{\Delta x}{2}$, $x_0 := 0$ and analogously for $y$. First of all, for $s \in \pm\frac12$, $r',s' \in \{ 0, \pm\frac12\}$ (edges),
\begin{align}
 B_{0,s}(x_{r'}, y_{s'}) \overset{\eqref{eq:tensorbasisedge}}{=} \frac23 B_{0}(x_{r'}) B_{s}(y_{s'}) = \frac23 B_0(0) \delta_{0,r'} \delta_{s,s'} = \delta_{0,r'} \delta_{s,s'}
\end{align}
since $B_0(0) = \frac32$. For the average, trivially $\int_C \tilde B_{r,s} \dd \vec x = 0$ since either $r$ or $s$ is $\pm\frac12$.

Second, for $(r,s) \in \left\{ \left(\pm\frac12,\pm\frac12\right) \right \}$ (nodes), and considering first $r',s' = \pm\frac12$
\begin{align}
 B_{r,s}(x_{r'}, y_{s'}) &\overset{\eqref{eq:tensorbasisnode}}{=} B_{r}(x_{r'}) B_{s}(y_{s'}) + \frac16 B_{r}(x_{r'}) B_{0}(y_{s'}) + \frac16 B_{0}(x_{r'}) B_{s}(y_{s'}) = \delta_{r,r'} \delta_{s,s'}
\end{align}
Next, if $s' = \pm\frac12$,
\begin{align}
 B_{r,s}(0, y_{s'}) &= B_{r}(0) \delta_{s,s'} + \frac16 B_{0}(0) \delta_{s,s'}= \left(-\frac14 + \frac16 \cdot \frac32 \right) \delta_{s,s'} = 0
\end{align}
since $B_{\pm\frac12}(0) = -\frac14$. A similar result is true for the perpendicular direction. We trivially have $\int_C \tilde B_{r,s} \dd \vec x = 0$ again. Finally, statement \eqref{eq:tensorbasisavg} is also clear.
\end{proof}

Consider now the tensor product of 1-d biorthogonal test functions, i.e. define
\begin{align}
 \tilde A_{r, s}(x,y) :=  A_r(x)  A_s(y) \qquad r,s \in \left\{\pm \frac12, 0 \right \}
\end{align}
with the normalization $\alpha = 0$ in \eqref{eq:centraltestfct1}--\eqref{eq:centraltestfct2}, i.e. with 
\begin{align}
 \int_{-\frac{\Delta x}{2}}^{\frac{\Delta x}{2}} A_{+\frac12} B_{+\frac12} \dd x  = \int_{-\frac{\Delta x}{2}}^{\frac{\Delta x}{2}} A_{-\frac12} B_{-\frac12} \dd x = \frac12
\end{align}
The test functions then are given by $A_0(x) = 1$ and \eqref{eq:1dbiorthtestsymm}.
As mentioned above, these tensor products are not the biorthogonal test functions in 2-d. However, it is now easily possible to compute the 2-d test functions $A_{rs}$. For example,

\begin{theorem} \label{thm:testfctsviatensor}
 The expansion of the 2-d test function $A_{+\frac12,0}$ in the basis of tensorial test functions $\tilde A$ is as follows:
 \begin{align}
  A_{+\frac12,0} &= \frac{8 \alpha_1 - 1 - \alpha_3}{2} \tilde A_{+\frac12,+\frac12}  + \frac{3(1+\alpha_3)}{2}\tilde A_{+\frac12,0} + \frac{8 \alpha_2 - 1 - \alpha_3}{2}  \tilde A_{+\frac12,-\frac12}  
 \end{align}
 in the notational conventions of the previous Section.
\end{theorem}
\begin{proof}
 $A_{+\frac12,0}$ cannot contain contributions from $\tilde A_{-\frac12,r}$, $r \in \left\{\pm \frac12, 0 \right \}$ since, on the one hand
 \begin{align}
  0= \int_C A_{+\frac12,0} B_{-\frac12,0} \dd \vec x = \frac23 \int_C A_{+\frac12,0} \tilde B_{-\frac12,0} \dd \vec x
 \end{align}
 and on the other, for $r = \pm\frac12$
 \begin{align}
  0&= \int_C A_{+\frac12,0} B_{-\frac12,r} \dd \vec x = \int_C A_{+\frac12,0} \tilde B_{-\frac12,r} \dd \vec x + \frac16 \int_C A_{+\frac12,0} \tilde B_{0,r} \dd \vec x + \frac16 \int_C A_{+\frac12,0} \tilde B_{-\frac12,0} \dd \vec x \\
  &= \int_C A_{+\frac12,0} \tilde B_{-\frac12,r} \dd \vec x + \frac16 \cdot \frac32 \int_C A_{+\frac12,0} B_{0,r} \dd \vec x + \frac16 \cdot 32 \int_C A_{+\frac12,0} B_{-\frac12,0} \dd \vec x \\
  &= \int_C A_{+\frac12,0} \tilde B_{-\frac12,r} \dd \vec x
 \end{align}
 $A_{+\frac12,0}$ cannot contain contributions from $\tilde A_{0,\pm\frac12}$ and $\tilde A_{0,0}$ either, for the same reasons. We are thus left with 
 \begin{align}
  A_{+\frac12,0} &= \gamma_1 \tilde A_{+\frac12,+\frac12} + \gamma_2 \tilde A_{+\frac12,0}+ \gamma_3 \tilde A_{+\frac12,-\frac12}
 \end{align}
for which
\begin{align}
  \alpha_1 &= \int_C A_{+\frac12,0} B_{+\frac12,+\frac12} \dd \vec x \\ &= \int_C A_{+\frac12,0} \tilde B_{+\frac12,+\frac12} \dd \vec x + \frac16 \int_C A_{+\frac12,0} \tilde B_{0,+\frac12} \dd \vec x + \frac 16 \int_C A_{+\frac12,0} \tilde B_{+\frac12,0} \dd \vec x \\&= \frac14 \gamma_1 + \frac1{12} \gamma_2\\
  \frac12 + \frac{\alpha_3}{2} &=\int_C A_{+\frac12,0} B_{+\frac12,0} \dd \vec x = \frac23  \int_C A_{+\frac12,0} \tilde B_{+\frac12,0} \dd \vec x = \frac13 \gamma_2\\
  \alpha_2 &= \int_C A_{+\frac12,0} B_{+\frac12,-\frac12} \dd \vec x \\
  &= \int_C A_{+\frac12,0} \tilde B_{+\frac12,-\frac12} \dd \vec x + \frac16 \int_C A_{+\frac12,0} \tilde B_{0,-\frac12} \dd \vec x +\frac16 \int_C A_{+\frac12,0} \tilde B_{+\frac12,0} \dd \vec x \\&= \frac14 \gamma_3 + \frac{1}{12} \gamma_2
 \end{align}
 which is a simple linear system for $\gamma_{1,2,3}$.
\end{proof}

\section{Numerical examples} \label{sec:numerical}

Next, some subtle differences in the results obtained with classical Active Flux discretizations and the new Petrov-Galerkin approach are illustrated. In all cases, RK3 with a CFL of 0.2 is used for the integration in time. 

\subsection{1-d Burgers' equation} \label{sec:numericalburgers}

\begin{figure}
 \centering
 \includegraphics[width=0.7\textwidth]{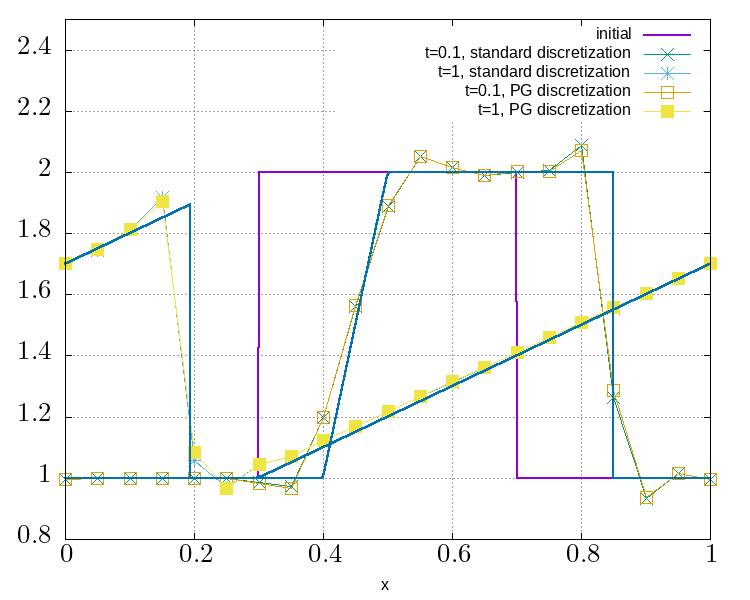}
 \caption{Comparison of results for classical upwind Jacobian splitting and the Petrov-Galerkin approach for Burgers' equation at times $t=0.1$ and $t=1$.}
 \label{fig:burgers}
\end{figure}

Consider 1-d Burgers' equation \eqref{eq:conslaw1d} with $f(q) = \frac{q^2}{2}$ and initial data
\begin{align}
 q_0(x) = \begin{cases} 2 & 0.4  < x < 0.6 \\
           1 & \text{else}
          \end{cases}
\end{align}
on a periodic domain $[0, 1]$ covered by a grid of 20 cells only. Figure \ref{fig:burgers} shows results obtained with the standard approach \eqref{eq:pointvaluesAF} and the Petrov-Galerkin discretization \eqref{eq:derivburgers}, with $\alpha_{i+\frac12} = 1$ since the data is uniformly positive. One observes only very subtle differences.

\subsection{Multi-d linear advection} \label{ssec:multidadvectiontest}

\begin{figure}
 \centering
 \includegraphics[width=0.32\textwidth]{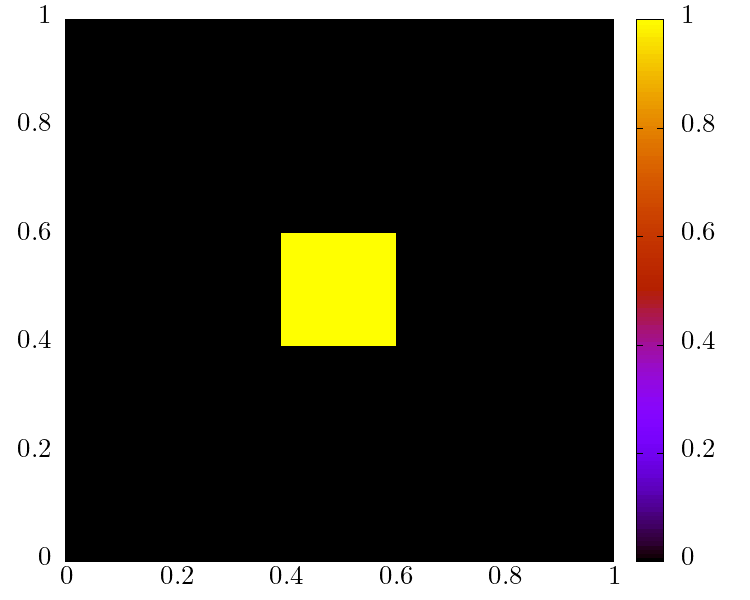} \hfill \includegraphics[width=0.32\textwidth]{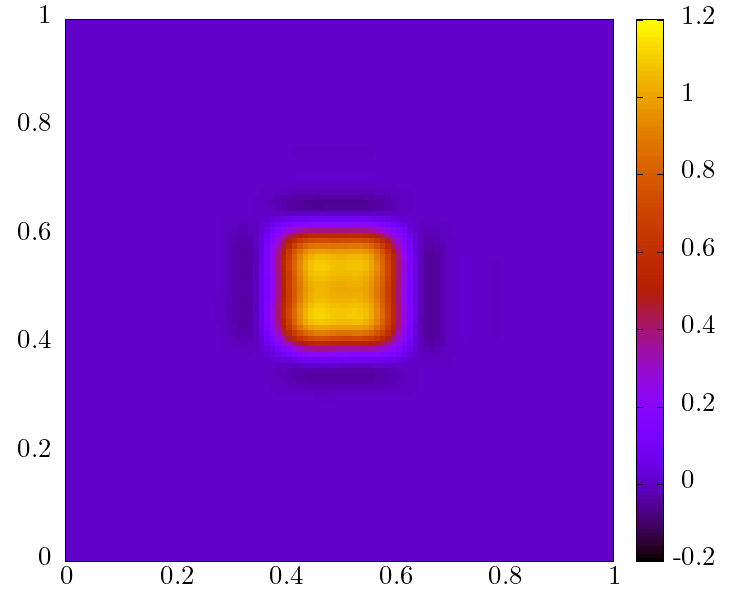} \hfill \includegraphics[width=0.32\textwidth]{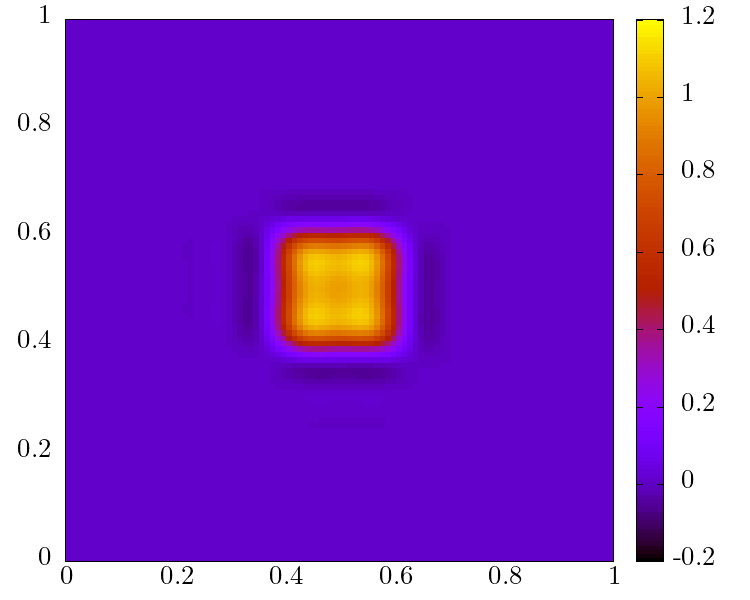} \\
 \includegraphics[width=0.7\textwidth]{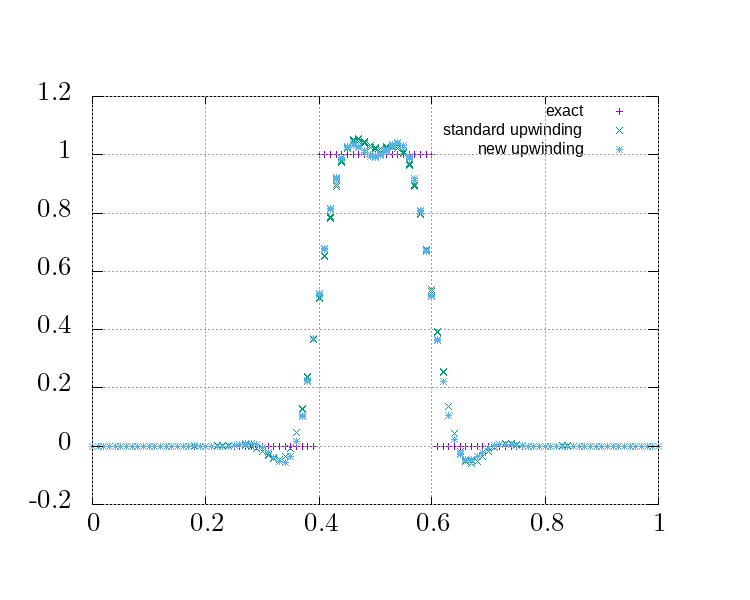} 
 \caption{Comparison of numerical results for linear advection with velocity $(1,\frac12)$ using standard upwinding and new upwinding, with parameter values as given in Section \ref{ssec:multidadvectiontest}. \emph{Top, from left to right}: Initial data / exact solution; solution at time $t=10$ using standard upwinding; solution at time $t=10$ using non-zero values in the new upwinding parameters. \emph{Bottom}: Slice through the simulation data at $y = 0.5$. One observes a slightly more symmetric result when making use of the new parameters.}
 \label{fig:square}
\end{figure}

Consider linear advection in two spatial dimensions
\begin{align}
 \del_t q + U^x \del_x q + U^y \del_y q = 0
\end{align}
with $U_x = 1$, $U^y = \frac12$, and initial data 
\begin{align}
 q_0(x, y) = \id_{[0.4, 0.6] \times [0.4, 0.6]}
\end{align}
on a domain $[0,1]^2$ with periodic boundaries, covered by a grid of $100 \times 100$ cells. Figure \ref{fig:square} shows the simulation result at time $t=10$ for the standard choice of upwinding, i.e.
\begin{align*}
 \alpha_3 &= 1 & \alpha_1 &= \alpha_2 = 0 &&\text{ in \eqref{eq:xderivexpand}}\\
 2(\alpha_{10} + \alpha_{11}) &= 1 & \alpha_1 &= \alpha_3 = \alpha_6 = \alpha_8 = 0 &&\text{ in \eqref{eq:nodederivx}} \\
 \alpha_9 &= 1 & \alpha_2 &= \alpha_4 = \alpha_5 = \alpha_7 = 0   &&\text{ in \eqref{eq:nodederivy}} 
\end{align*}
and for the choice
\begin{align*}
 \alpha_3 &= 1 &  \alpha_2 &= 0.2 & \alpha_1 &= 0 &&\text{ in \eqref{eq:xderivexpand}}\\
 2(\alpha_{10} + \alpha_{11}) &= 1 & \alpha_1 &= \alpha_3 = 0.2 & \alpha_6 &= \alpha_8 = 0&&\text{ in \eqref{eq:nodederivx}} \\
 \alpha_9 &= 1 & \alpha_2 &= \alpha_4 = 0.2 & \alpha_5 &= \alpha_7 = 0  &&\text{ in \eqref{eq:nodederivy}} 
\end{align*}
Here, the additional jumps terms are included in upwind direction only. The value 0.2 is not based on any fundamental principle. The experimental results seem to indicate a possibility to obtain more symmetric results upon making use of the new parameters. A detailed study of their influence is, however, subject of future work. 

\section{Conclusions and outlook}

So far, Active Flux was employing a continuous approximation known from continuous Finite Elements, an average update similar to Finite Volume/Discontinuous Galerkin methods, and (in the semi-discrete case) a finite-difference-type update of the point values. An obvious difference to both CG and DG is the absence of (even a local) mass matrix in Active Flux. This work shows that the update equations of the semi-discrete Active Flux method can be obtained upon using a biorthogonal set of test functions in a Petrov-Galerkin framework. The cases studied here explicitly are arbitrarily high-order accurate Active Flux methods in one spatial dimension and the classical third-order accurate Active Flux method on two-dimensional Cartesian grids. This work also can be seen as an example of beneficial use of biorthogonal test functions, which by construction give a method without a mass matrix and 
through their discontinuities reproduce upwind derivative choices for stabilization.
In fact, the new framework even suggests additional jump terms, which so far have not been considered in AF literature but might be useful in future. With globally continuous basis functions and discontinuous test functions, AF thus indeed seems to lie on the boundary between CG and DG. This work thus contributes to placing the rather recent Active Flux amidst well-known established methods and shows its relation to them.

Future work will be devoted to a detailed study how the additional jump terms that have appeared in the presented Petrov-Galerkin framework can be used for stabilization, to a more detailed study of the multi-dimensional case, covering in particular arbitrary orders of accuracy, and to an extension of this formulation to other types of meshes. Triangular meshes seem particularly promising, in view of studies such as \cite{abgrall25dg}, which appeared after this work, and since biorthogonal polynomial systems for triangles have received considerable attention in the past.

\section*{Acknowledgement}

% #############################################################################################################################
The author thanks Davide Torlo and Mario Ricchiuto for many helpful discussions.
% #############################################################################################################################

\bibliographystyle{alpha}
\newcommand{\etalchar}[1]{$^{#1}$}

\end{document}